\documentclass[12pt]{article}

\usepackage{amsmath}
\usepackage{amssymb}

\usepackage[T1]{fontenc} 
\usepackage{sets pace}
\usepackage{comment}
\usepackage{colortbl}
\definecolor{text1}{cmyk}{1,.65,0,0} % blue text colour
\definecolor{text2}{rgb}{1,0,0} % red text colour
\definecolor{text3}{cmyk}{0,0,0,1} % black text colour
\definecolor{text4}{cmyk}{0,0,0,0.5} % grey text colour
\definecolor{text5}{cmyk}{1.0,0.0,1.0,0} % green text colour

\usepackage{amsfonts}
\usepackage{amsthm}
\usepackage{graphicx}
\usepackage{caption}
\usepackage{subcaption}
\usepackage{latexsym}
\usepackage[most]{tcolorbox}
\makeatletter
\renewcommand{\@seccntformat}[1]
{\csname the#1\endcsname.\enspace}
\makeatother
\newtheorem{definition}{Definition}
\newtheorem{theorem}{Theorem}
\newtheorem{lemma}{Lemma}
\newtheorem{remark}{Remark}
\newtheorem{corollary}{Corollary}
\newtheorem{example}{Example}

\begin{document}
%\doublespacing
\begin{center}
   {\bf Bayesian inference and prediction for mean-mixtures of normal
distributions \footnote{\today} }
\end{center}

\begin{center}
{\sc Pankaj Bhagwat$^{a}$ \& \'{E}ric Marchand
} \\

{\it a  Universit\'e de Sherbrooke, D\'epartement de math\'ematiques, Sherbrooke Qc,    CANADA, J1K 2R1 \quad 
(e-mails: pankaj.uttam.bhagwat@usherbrooke.ca; eric.marchand@usherbrooke.ca)} \\

\end{center}

\vspace*{0.2cm}
\begin{center}
{\sc Abstract} \\
\end{center}
\vspace*{0.1cm}
\small  
We study frequentist risk properties of predictive density estimators for mean mixtures of multivariate normal distributions, involving an unknown location parameter $\theta \in \mathbb{R}^d$, and which include  multivariate skew normal distributions.   We provide explicit representations for Bayesian posterior and predictive densities, including  the benchmark minimum risk equivariant (MRE) density, which is minimax and generalized Bayes with respect to an improper uniform density for $\theta$.   For four dimensions or more, we obtain Bayesian densities that improve uniformly on the MRE density under Kullback-Leibler loss.   
We also provide plug-in type improvements, investigate implications for certain type of parametric restrictions on $\theta$, and illustrate and comment the findings based on numerical evaluations.   

\vspace*{0.2cm}


\normalsize

\section{Introduction}
The findings of this paper relate to predictive density estimation for mean mixture of normal distributions.    The modelling of data via mixing multivariate normal distributions has found many applications and lead to methodological challenges for statistical inference.    These include finite mixtures, as well as continuous mixing on the mean and/or the variance.   Whereas, scale or variance mixtures of multivariate normal distributions compose a quite interesting subclass of spherically symmetric distributions, modelling asymmetry requires mixing on the mean and prominent examples are generated via variance-mean mixtures (e.g., \cite{barndorff-nielsenetal.1982}), as well as mean mixtures of multivariate normal distributions (e.g., \cite{ambj2021, arvalle+azzalini2021}) and references therein).      Moreover, such mean mixtures, which are the subject of study here, generate or are connected to multivariate skew-normal distributions (e.g., \cite{azz+dalla1996}) which have garnered much interest over the years.

The development of shrinkage estimation techniques, namely since Stein's inadmissibility finding (\cite{stein1956}) concerning the maximum likelihood or best location equivariant estimator under squared error loss in three dimensions or more, has had a profound impact on statistical theory, thinking, methods, and practice.  Examples include developments on sparsity and regularization methods, empirical Bayes estimation, adpative inference, small area estimation, non-parametric function estimation, and predictive density estimation.   Cast in a decision-theoretic framework, Stein's original result has been expanded in many diverse ways, namely to other distributions or probability models, and namely for spherically symmetric and elliptically symmetric distributions (see for instance, \cite{fsw2018}).  There have been fewer findings for multivariate skew-normal or mean mixtures of normal distributions, but the recent work of Kubokawa et al. \cite{ksy2020} establishes point estimation minimax improvements of the best location equivariant estimator under quadratic loss, when the dimension of the location parameter is greater than or equal to four, and with underlying known perturbation parameter. 

Predictive density estimation has garnered much interest over the past twenty years or so, and addresses fundamental issues in statistical predictive analysis.  Decision-theoretic links between shrinkage point estimation and shirnkage predictive density estimation for normal models have surfaced (e.g., \cite{komaki2001}, \cite{glx2006}) and stimulated much activity (see for instance \cite{banff}), including findings for restricted parameter spaces (e.g., \cite{fmrs2011}).  
The main objective of this work is thus to explore the problem of predictive density estimation for mean mixtures of normal (MMN) distributions.  A secondary objective is to provide novel representations for Bayesian posterior distributions and predictive densities, which have been found to be lacking in the literature.  

Following early findings of Komaki (e.g., \cite{komaki2001}) on the predictive density estimation problem for multivariate normal models under Kullback-Leibler loss, George, Liang and Xu in \cite{glx2006} exhibited further parallels with the point estimation problem for normal distribution under quadratic loss. They provide sufficient conditions on marginal distributions and prior distributions to get improved shrinkage predictive density estimators when the dimension is greater than or equal to three. Thus, motivated by these connections, it is interesting to investigate whether such shrinkage plays any role in the predictive density estimation problem for mean-mixture of multivariate normal models and we focus on frequentist risk efficiency of predictive density estimators under Kullback-Leibler loss.   Our contribution here consists in identifying classes of plug-in type predictive densities and  of Bayes predictive densities which are minimax and dominate the benchmark minimum equivariant estimator (MRE) for the case when the dimension of the location parameter is greater than or equal to four.

The organization of this manuscript is as follows.   Section 2.1 contains several introductory definitions, properties and examples of MMN models, including a useful canonical form which subdivides the MMN distributed vector into $d$ independent components, one of which  a univariate MMN distribution and the others as normal distributions.   Section 2.2 focuses on the predictive estimation framework with a KL loss decomposition, and an initial representation for the MRE density accompanied by various examples.  Section 2.3 expands on the calculation of minimax risk and a representation in terms of the entropy of a univariate distribution.   Section 3 is devoted to Bayesian posterior and predictive analysis with several novel representations.   Sections 4.1 and 4.2,  namely Theorems \ref{theoremdominanceplugins}, Theorem \ref{theoremmaindominance} and Corollary \ref{superharmonic conditions}, contain the main dominance findings, with plug-in type and Bayesian improvements.    In both cases, the main technique employed rests upon the canonical transformation presented in Section 2.1 and permits to split up the KL risk as the addition of two parts, one of which can be operated on using known normal model prediction analysis findings.    Section 4.3 deals with parametric restrictions and further applications of Theorems \ref{theoremdominanceplugins} and \ref{theoremmaindominance}.    Finally, some further  illustrations are provided in Section 5.

\section{Preliminary results and definitions }

Here are some details, properties and definitions on mean mixture of normal distributions, its canonical form, and predictive density estimation.  In the following, we will denote $\phi_d(z;\Sigma)$ the probability density function (pdf) of a $N_d(0,\Sigma)$ distribution evaluated at $z \in \mathbb{R}^d$ and for positive definite $\Sigma$.   When $\Sigma=I_d$, we may simplify the writing to $\phi_d(z)$, and then for $d=1$ to $\phi(z)$.    We will denote $\Phi$ the cdf of a $N(0,1)$ distribution.   

\subsection{The model}

The distributions of interest are mean-mixtures of multivariate normal distributions, both for our observables and densities to be estimated by a predictive density estimator.   Such distributions connect to multivariate skew-normal distributions and have been the object of interest in recent work with studies of stochastic properties (e.g., \cite{ambj2021}, \cite{arvalle+azzalini2021}), and shrinkage estimation about its location parameter (\cite{ksy2020}).

\begin{definition}
\label{defmmn}
A random vector $X \in \mathbb{R}^d$ is said to have a mean-mixture of normal distributions (MMN), denoted as $X \sim MMN_d(\theta,a, \Sigma , \mathcal{L})$, if it admits the representation
\begin{equation}
    \label{model} X|V=v \sim N_d(\theta + v a\,, \Sigma) \,, \,     V  \sim \mathcal{L}\,,
 \end{equation}
where $\theta \in \mathbb{R}^d $ is a location parameter, $a \in \mathbb{R}^d - \{0\}$ is a known perturbation vector, $\Sigma$ is a known positive definite covariance matrix, and $V$ is a scalar random variable with cdf $\mathcal{L}$.
\end{definition}

Alternatively, the random vector $X$ has stochastic representation 
\begin{align}
    X = \theta + \Sigma^{1/2}Z + V a \,,
\end{align}
where $Z \sim N_d(0,I_d)$ and $V \sim \mathcal{L}$ on $\mathbb{R}$, and its 
probability density function can be expressed as:
\begin{align}
    \nonumber p(x|\theta) &=  \mathbb{E}^V \left\lbrace\phi_d\left(x - \theta - Va,  \Sigma\right) \right\rbrace \\
    \label{mmndensityequiv.form}   & =   \phi_d\left(x - \theta,  \Sigma\right) \, \mathbb{E}^V \left( e^{-\frac{V^2}{2}a^T\Sigma^{-1}a}e^{V \, (x-\theta)^T\Sigma^{-1}a } \right) \,.
\end{align}
Thus, we note that the density function of MMN random vector can be decomposed in two parts: one symmetrical density $\phi_d(\cdot)$ and the other part which is a function of the projection of $(x - \theta)$ in the direction of $\Sigma^{-1}a$.   Moreover, this construction isolates the asymmetry in the direction $\Sigma^{-1}a$ and the scale is controlled by the random variable $V$. 

\begin{remark}
\label{closure}
It is easy to see that the family of MMN distributions is closed under linear combinations of independent components.   Specifically, if $X_i|\theta \sim MMN_d(\theta,a, \Sigma_i , \mathcal{L}_i)$, $i=1,\ldots, n$, are independently distributed, then $\sum_{i=1}^n b_i X_i|\theta \sim MMN_d((\sum_{i=1}^n b_i) \, \theta, a, \sum_{i=1}^n b_i^2 \, \Sigma_i, \mathcal{L}_0)$ with $\mathcal{L}_0$ the cdf of the mixing variable $V_0 =^d \sum_{i=1}^n b_i V_i$.  Namely, for the identically distributed case with $\Sigma_i=\Sigma$ and the sample mean with $b_i=1/n$, we obtain that  
\begin{equation}
\nonumber
\bar{X}|\theta \, \sim \, MMN_d(\theta, a, \Sigma/n, \mathcal{L}_0) \,, \hbox{ with }
\mathcal{L}_0\, \hbox{ the cdf of } \bar{V}\,.
\end{equation}

\noindent  It thus follows, as observed in \cite{ksy2020}, that findings applicable for a single MMN distributed observable $X$ can be extended to the random sample case.   
\end{remark}

We now turn our attention to a fundamental decomposition, or canonical form, (\cite{ambj2021}) for MMN distributions which will be most useful.   

\begin{lemma}
\label{lemmacanonical}
For a random vector $X \sim MMN_d(\theta,a,\Sigma, \mathcal{L})$ as in (\ref{model}), there exists an orthogonal matrix $H$ such that the first row of $H$ is proportional to $a^{\top}\, \Sigma^{-1/2}$ and $Z = H\Sigma^{-1/2}X $ has a $MMN_{d}(H\Sigma^{-1/2}\theta,a_0, I_d, \mathcal{L})$ distribution with $a_0 = (\sqrt{a^T\Sigma^{-1}a}, 0,\dots,0)^T$. 
\end{lemma}

Such a $Z$ may be referred to as a canonical form and is comprised of $d$ independent components. Moreover $Z-H\Sigma^{-1/2}\theta$ has $d-1$ components which are $N(0,1)$ distributed and another distributed as $MMN_1(0,a_0,1,\mathcal{L})\,$.   Such a canonical form construction is not unique and depends on the choice of $H$. 

\noindent 
As already mentioned, the family of MMN distributions contains many interesting distributions and we refer to the above-mentioned references for various properties.   We expand here with illustrations, which will also inform us for our predictive density problem and related Bayesian posterior analysis.   A prominent example is the multivariate skew normal distribution due to Azzalini and Dalla Valle \cite{azz+dalla1996}. If we consider $V \sim TRN(0,1),$ the standard truncated normal distribution on $R_{+}$ in (\ref{model}), we get the multivariate skew-normal family of distributions with densities
  \begin{equation}
  \label{skewnormaldensity}
    p(x|\theta)  
     \, = \, 2\phi_d\left(x - \theta; \Sigma + aa^T \right)\Phi\left( \frac{(x - \theta)^{\top}\Sigma^{-1}a}{\sqrt{1 + a^{\top}\Sigma^{-1}a}} \right).
\end{equation}
We denote this as $X \sim SN_d(\theta,a,\Sigma)$. Here, we note that $V \sim \sqrt{\chi^2_1}$, i.e. the square root of a Chi-square  distribution with $k=1$ degrees of freedom.  Various other choices of the mixing density have appeared in the literature (e.g., \cite{arvalle+azzalini2021}), namely cases where $V \sim \sqrt{\chi^2_k}$ or $V$ is Gamma distributed.  Here is a general result containing such cases as well as many others. 

\begin{theorem}
\label{theorem-mixing}
For a mixing density of the form
\begin{align}
    \ell(v) = h(v)\,e^{-vc_2}\,e^{-\frac{v^2}{2}c_1} \, \mathbb{I}_{(0,\infty)}(v) \,,
    \label{mixinggeneral}
\end{align}
with $c_1 > 0,c_2 \in \mathbb{R}$ or $c_1=0,c_2 \geq 0$, the corresponding pdf of $X$ in (\ref{model}) is given by 
\begin{align}
    \label{densitygeneral} p(x|\theta) & =  \frac{1}{c_1'} \, \phi_d\left(x -\theta,  \Sigma\right) \, \frac{\mathbb{E}\left[\left.h\left\lbrace \frac{1}{c_1'}\left( Z + \frac{c_2'}{c_1'}\right)\right\rbrace \right\vert Z + \frac{c_2'}{c_1' }\geq 0\right] }{R\left(\frac{c_2'}{c_1' }\right)},
\end{align}
 with $Z \sim N(0,1)$, $c_1' = \left( c_1 + a^{\top}\Sigma^{-1}a\right)^{1/2} $, $c_2' = (x - \theta)^{\top}\Sigma^{-1}a - c_2 \,$, and $R(\cdot)$ the reverse Mill's ratio given by $R(t)=\phi(t)/\Phi(t), t \in \mathbb{R}$. 
%with  $\frac{c_2'}{\sqrt{c_1'}} = \frac{(x - \theta)^{\top} \, \Sigma^{-1}a - c_2}{\sqrt{c_1 + a^{\top}\Sigma^{-1}a}}$ and $c_1' = c_1 + a^{\top}\Sigma^{-1}a$. 

{\bf Proof}.  The result follows from (\ref{mmndensityequiv.form}) as 
\begin{eqnarray}
    \nonumber \mathbb{E}^V \left( e^{-\frac{V^2}{2}a^T\Sigma^{-1}a}e^{V(x-\mu)^T\Sigma^{-1}a}\right) & = & \int\limits_{0}^{\infty} \, e^{-\frac{v^2}{2} (c_1')^2}e^{vc_2'}\,h(v)\,dv \\
    \nonumber & = & \frac{\sqrt{2\pi}}{c_1'} \, \, e^{\frac{c_2'^2}{2c_1'^2}} \int\limits_{0}^{\infty} h(v) \, \frac{c_1'}{\sqrt{2\pi}} \; e^{-\frac{c_1'^2 \left( v - \frac{c_2'}{(c_1')^2}\right)^2}{2}}\, dv \\
    \nonumber & = & \frac{\sqrt{2\pi}}{c_1'} \; e^{\frac{c_2'^2}{2c_1'^2}} \; \mathbb{E}\left\lbrace\left.h\left(\frac{Z}{c_1'} + \frac{c_2'}{c_1'^2}\right) \right\vert Z + \frac{c_2'}{c_1'}\geq 0\right\rbrace \, \Phi\left(\frac{c_2'}{c_1' }\right)\,.  \qed
\end{eqnarray}
\end{theorem}

We point out that the above Theorem applies for $c_1=c_2=0$ and thus covers all absolutely continuous distributions on $\mathbb{R}_+$.   Here are nevertheless specific examples of Theorem \ref{theorem-mixing} and model density (\ref{densitygeneral}).

\begin{example}
\label{examplemixing}
\begin{enumerate}
\item[ {\bf (A)}]   Gamma mixing with $h(v) \, = \,  \frac{v^{\alpha - 1}}{\Gamma(\alpha) \beta^{\alpha}}$.   Theorem \ref{theorem-mixing} applies with $c_1=0$ and $c_2=1/\beta$, and the model density is given by (\ref{densitygeneral}) with the above $h$, $c_1'\,=\, \left( a^{\top} \Sigma^{-1} a \right)^{1/2}$ and $c_2'\,=\, (x - \theta)^{\top}\Sigma^{-1}a  - ( 1/\beta)$.   The density was studied in \cite{ambj2021,adcock+shutes}.  The exponential case with $\alpha=1$ simplifies with 
\begin{equation}
\noindent
    p(x|\theta) = \frac{1}{\beta c_1'} \,\, \frac{\phi_d\left(x - \theta; \Sigma  \right)}{ R\left(\frac{c_2'}{c_1'}\right)}\,.
    \label{expmixing}
\end{equation}
More generally for positive integer $\alpha$, the density's expression brings into play the $(\alpha-1)^{\hbox{th}}$ lower-truncated moment of a normal distribution.    For instance, with $\mathbb{E}\left\lbrace(Z+\Delta)|Z+\Delta \geq 0\right\rbrace \, = \, \Delta \, + \, R(\Delta),$ we obtain for the case $\alpha=2$ the model density:

\begin{equation}
\nonumber
p(x|\theta) \, = \, \frac{\phi_d(x-\theta, \Sigma)}{(c_1' \, \beta)^2} \, 
\left\lbrace   \frac{c_2'/c_1'}{R(c_2'/c_1')} \, + \, 1\right\rbrace \,,
\end{equation}
with the above $c_1'$ and $c_2'$.

\item[ {\bf (B)}]   $\sqrt{\chi_k^2}$ mixing with $h(v) \, = \, \frac{(\frac{1}{2})^{k/2 -1}}{\Gamma(k/2)} \, v^{k-1}$, $c_1=1$, $c_2=0$, and $k>0$.   The corresponding model density is given by (\ref{densitygeneral}) with the above $h$, $c_1' = \left(1+a^{\top} \Sigma^{-1} a   \right)^{1/2}$, and $c_2' \, = \, (x - \theta)^{\top}\Sigma^{-1}a$.

%$$p(x|\theta) \, = \, \frac{2\sqrt{\pi}}{\Gamma(k/2)}\frac{\phi_d\left(x - \theta; \Sigma + aa^T \right)}{\left(\sqrt{2(1 + a^{\top}\Sigma^{-1}a)}\right)^{k-1}}\Phi\left( \frac{(x - \theta)^{\top}\Sigma^{-1}a}{\sqrt{1 + a^{\top}\Sigma^{-1}a}} \right)\\ \times \mathbb{E}\left.\left[\left(Z + \frac{(x - \theta)^{\top}\Sigma^{-1}a}{\sqrt{1 + a^{\top}\Sigma^{-1}a}}\right)^{k-1}\right|Z + \frac{(x - \theta)^{\top}\Sigma^{-1}a}{\sqrt{1 + a^{\top}\Sigma^{-1}a}} \geq 0 \right].$$

The density was given in \cite{arvalle+azzalini2021} and, as previously noted, the case $k=1$ reduces to the skew-normal case in (\ref{skewnormaldensity}).   As in Example {\bf (A)} for positive integer $k$, the density's expression involves a lower-truncated moment of a normal distribution.
 
\item[ {\bf (C)}]  Kummer type II mixing with $c_2=c/\sigma$, $c_1=0$, $h(v) = \frac{\sigma^b}{\Gamma(a) \, \psi(a,1-b,c)} \, \frac{v^{a-1}}{(v+\sigma)^{a+b}}$ with $a, c, \sigma >0$, $b \in \mathbb{R}$, and $\psi$ the confluent hypergeometric function of type II defined for $\gamma_1, \gamma_3>0$ and $\gamma_2 \in \mathbb{R}$ as $\psi(\gamma_1, \gamma_2, \gamma_3) \, = \, \frac{1}{\Gamma(\gamma_1)} \, \int_{\mathbb{R}_+} t^{\gamma_1-1} (1+t)^{\gamma_2-\gamma_1-1} \, e^{-\gamma_3 t} \, dt$.   This class of densities includes for $b=-a$ the Gamma densities in {\bf (A)}, as well as Beta type II densities for $c=0$ and $b>0$  The resulting mean-mixture density is given by (\ref{densitygeneral}) and involves interesting expectations of the form $\mathbb{E} \left( \frac{W^{a-1}}{(W+\sigma)^{a+b}} \vert W \geq 0 \right)$ where $W \sim N(\Delta, 1)$ with $\Delta=c_2'/c_1'$.  
\end{enumerate}
\end{example}

\subsection{The prediction problem}
\label{PredictionProblem}

Consider $X|\theta \sim MMN_d(\theta,a, \Sigma_X , \mathcal{L}_1)$ and $Y|\theta \sim MMN_d(\theta,a, \Sigma_Y , \mathcal{L}_2),$ independently distributed as in Definition \ref{defmmn}, i.e.
\begin{equation}
\label{modelXY}    X|\theta,V_1 \sim N_d(\theta + V_1\;a, \Sigma_X)\,, \,  
  Y|\theta,V_2 \sim N_d(\theta + V_2\;a, \Sigma_Y)\,, \hbox{ with } V_1 \sim \mathcal{L}_{1}\,,\, V_2 \sim \mathcal{L}_{2}\,.  
\end{equation}
Let $p(x|\theta)$ and $q(y|\theta)$ denote the conditional densities of $X$ and $Y$ given $\theta$, respectively.
Based on observing $X=x$, we consider the problem of finding a suitable predictive density estimator $\hat{q}(y;x)$ for $q(y|\theta)\,, y \in \mathbb{R}^d\,.$

The ubiquitous Kullack-Leibler (KL) divergence between two Lebesgue densities $f$ and $g$ on $\mathbb{R}^m$, defined as 
\begin{equation}
\nonumber  \rho(f,g) \, = \, \int_{\mathbb{R}^m} f(t) \, \log \frac{f(t)}{g(t)} \, dt\,,
\end{equation}
is the basis of Kullback-Leibler loss given by 
\begin{equation}
\label{klloss}  L(\theta, \hat{q}) \, = \, \rho(q_{\theta}, \hat{q})\,.
\end{equation}
We will make use of Lemma \ref{lemmacanonical}'s canonical form as in (\ref{lemmacanonical}) to transform a mean mixture of normal distributions vector into two independent components and to capitalize on the corresponding simplification for KL divergence which is as follows.

\begin{lemma}
\label{lemmaKL}
Let  $T=(T_{(1)}, T_{(2)}) \in \mathbb{R}^m$ and $U=(U_{(1)}, U_{(2)}) \in \mathbb{R}^m$ be random vectors subdivided into independently distributed components $T_{(i)}$ and $U_{(i)}$ of dimensions $m_i$ for $i=1,2$ with $m_1+m_2=m$.   Denote $f$ and $g$ the densities of $T$ and $U$, respectively, and $f_1, f_2, g_1, g_2$ the densities of $T_{(1)}, T_{(2)}, U_{(1)}, U_{(2)}$, respectively.  Then, we have
\begin{equation}
\label{rhodecomposition}   \rho(f,g) \, = \, \rho(f_1,g_1) \, +  \, \rho(f_2,g_2)\,.
\end{equation}
{\bf Proof.}  By independence, we have 
\begin{equation}
\nonumber  \rho(f,g) \, = \, \mathbb{E}^T \left\lbrace   \log \left(\frac{f_1(T_1) f_2(T_2)}{g_1(T_1) g_2(T_2)} \right) \right\rbrace \, = \, \mathbb{E}^T \left\lbrace   \log \left(\frac{f_1(T_1) }{g_1(T_1) } \right) \right\rbrace \, + \, \mathbb{E}^T \left\lbrace   \log \left(\frac{f_2(T_2) }{g_2(T_2) } \right) \right\rbrace \,  ,
\end{equation}
which is (\ref{rhodecomposition}).  \qed 
\end{lemma}
 We evaluate the performance of the density estimators using KL loss in (\ref{klloss}),

and the associated KL risk function
\begin{align}
    R_{KL}(\theta, \hat{q}) = \int_{\mathbb{R}^d} \{\int_{\mathbb{R}^d} q(y|\theta)\;\text{log}\;\frac{q(y|\theta)}{\hat{q}(y;x)}\; dy \}\; p(x|\theta)\; dx.
\end{align}

For a prior density $\pi$ for $\theta$ with respect to a $\sigma-$finite measure $\nu$, it is known (e.g.,  \cite{ad1975,aitchison1975}) that the Bayes predictive density is given by
\begin{align}
    \hat{q}_{\pi}(y;x) = \int_{\mathbb{R}^d} q(y|\theta)\; p(x|\theta)\;\pi(\theta)\;d\nu(\theta).
    \label{eq:Bayespredictivedensityfor pi}
\end{align}

A benchmark predictive density estimator for $q(y|\theta), y \in \mathbb{R}^d$, is given by the Bayes predictive density estimator $\hat{q}_U(y; X), y \in \mathbb{R}^d$, with respect to the uniform prior density on $\mathbb{R}^d$.  It is known to be the minimum risk equivariant (MRE) predictive density estimator under changes of location, as well as minimax.    In \cite{kms2015}, a representation, which applies to both integrated squared-error loss and KL loss, for $\hat{q}_U$ is provided.   For our prediction problem,  the following result makes use of this representation and  summarizes the above optimality properties.  

\begin{lemma}
\label{theoremmre}
 The MRE predictive density estimator of the density of $Y$ relative to model (\ref{modelXY}) under KL loss, is  given by the Bayes predictive density $\hat{q}_U$ under prior $\pi_U(\theta) =1\,$.   Furthermore, we have

\begin{align}
    \hat{q}_{U}(\cdot;X) \sim MMN_d(X,a, \Sigma_X + \Sigma_Y , \mathcal{L}_3)\,,
    \label{mreestimator}
\end{align}
where $\mathcal{L}_3$ is the cdf of $V_3=V_2-V_1$.
Finally, $\hat{q}_{U}(y;X)$ is minimax under KL loss.
\end{lemma}
{\bf Proof.}  The MRE and minimax properties are given in \cite{murray1977} and \cite{lb2004}, respectively.  For a location family prediction problem with $X \sim p(x-\theta)$ and  $Y \sim q(y- \theta)$ independently distributed, it is shown in \cite{kms2015} that 
\begin{equation}
\nonumber   \hat{q}_{U}(y;X) \, = \, q * \bar{p}(y-x) \,,\hbox{ with } \bar{p}(t)=p(-t)\,,
\end{equation} 
i.e., the convolution of $q$ and the additive inverse of $p$ followed by a change of location equal to $x$.   For model (\ref{model}), the above convolution $q * \bar{p}$ is given by the density of $Y-X$ in model (\ref{model}) with $\theta=0$, and the result follows since
\begin{equation}
\nonumber
Y-X|V_1, V_2 \sim N_d((V_2-V_1)\, a, \Sigma_X+ \Sigma_Y)\,. \;\;\;\;\;\;\;\;\;\;\;\;\;\;\;\;\;\; \qed   
\end{equation}

Here, we can see that the MRE density estimator also belongs to the class of MMN distributions with same perturbation parameter $a$ and location parameter $x$.    
As well, the distribution of the difference $V_2-V_1$ plays a key role in Theorem \ref{theoremmre}'s representation of the MRE predictive density, and as illustrated in the next subsection of examples.

\subsection{Minimax risk and entropy}

The Kullback-Leibler risk expressions brings into play the entropy associated with MMN distributions.   Such a measure is not easily manipulated into a closed form (see for instance \cite{entropy2012} for the study of entropy for skewed-normal distributions), but they can be expressed in terms of the entropy of a univariate MMN distribution, as illustrated with the following expansion of the constant and minimax risk of the MRE density $\hat{q}_U$ in the context of model (\ref{modelXY}).    For a Lebesgue density on $\mathbb{R}^d$, defined as 
\begin{equation}
\nonumber   H(f) \, = \, - \int_{\mathbb{R}^d} f(t) \, \log f(t) \, dt\,,
\end{equation}
we will make use of the following well-known and easily established properties.

\begin{lemma}
\label{lemmaentropy}
\begin{enumerate}
\item[ {\bf (a)}]  For $T \in \mathbb{R}^d$ with density $f$ and $U=\psi(T) \sim g$ with $\psi: \mathbb{R}^d \to \mathbb{R}^d$ invertible with inverse Jacobian $J_{\psi}$, we have 
 $  H(g) \, = \, - \mathbb{E} \log |J_{\psi}| + H(f)\,; $ 
\item[ {\bf (b)}]  Let $T=(T_{(1)},T_{(2)}) \sim f$ be a random vector with independently distributed components $T_{(1)} \sim f_1$ on $\mathbb{R}^{m_1}$ and $T_{(2)} \sim f_2$ on $\mathbb{R}^{m_2}$.  Then (as in Lemma \ref{lemmaKL}), we have $H(f) \, = \, H(f_1) \,+\, H(f_2)\,$. 
\end{enumerate}

\end{lemma}

As implied by part {\bf (a)} of the above lemma, the entropy $H(f_{\mu})$ is constant as a function of $\mu$ for location family densities $f_{\mu}(t) = f_0(t- \mu)$, as is the case for 
$MMN_d(\mu, b, \Sigma, \mathcal{L})$ densities.  Now, we have the following dimension reduction decomposition for the entropy $H_d(b, \Sigma, \mathcal{L})$ of a $MMN_d(0, b, \Sigma, \mathcal{L})$ density.

\begin{lemma}
\label{lemmaentropymmn}
We have for $d \geq 2$: $$H_d(b,\Sigma, \mathcal{L}) \,  = \, H_1(\sqrt{b^{\top} \Sigma^{-1} b}, 1, \mathcal{L}) \, + \, \frac{d-1}{2} \, \{1 + \log (2\pi)\}\, + \, \frac{1}{2} \log |\Sigma|\,.$$
\end{lemma}
{\bf Proof.}   Let $X \sim MMN_d(0, b, \Sigma, \mathcal{L}) $, which has entropy $H_d(b, \Sigma, \mathcal{L})$, and set $Z=H \, \Sigma^{-1/2} \, X \sim f_Z$ with $H$ orthogonal having first row $\frac{b^{\top} \Sigma^{-1/2}}{\sqrt{(b^{\top} \Sigma^{-1} b)}}$.   It follows from part {\bf (a)} of Lemma \ref{lemmaentropy} that $H(f_Z) = - \frac{1}{2} \, \log |\Sigma| + H_d(b, \Sigma, \mathcal{L})$.  
From Lemma \ref{lemmacanonical}, we have $Z=(Z_1, Z_{(2)})^{\top}$ with $Z_1 \sim MMN_1(0,\sqrt{(b^{\top} \Sigma^{-1} b)}, 1, \mathcal{L})$ and $Z_2 \sim N_{d-1}(0,I_{d-1})$ independently distributed, and the result follows from part {\bf (b)} of Lemma \ref{lemmaentropy} and a straightforward evaluation of the entropy $H(\phi_{d-1})$. \qed

With the above, we conclude with an expression for the constant and minimax risk.

\begin{theorem}  In the context of model (\ref{modelXY}), the Kullback-Leibler risk of the MRE density $\hat{q}_U$ is given by
\begin{equation}
\label{riskmre}
R_{KL}(\theta, \hat{q}_U) \, = \, H_1(\sqrt{a^{\top} \Sigma_S^{-1} a}, 1, \mathcal{L}_3) \, - \, H_1(\sqrt{a^{\top} \Sigma_Y^{-1} a}, 1, \mathcal{L}_2) + \frac{1}{2} \log \frac{\Sigma_S}{\Sigma_Y}\,,
\end{equation}
with $\Sigma_S \, = \, \Sigma_X + \Sigma_Y$.
\label{theoremminimaxrisk}
\end{theorem}
{\bf Proof.}
We have for $\theta \in \mathbb{R}^d\,$
\begin{eqnarray*} 
R_{KL}(\theta, \hat{q}_U) \, & = & \, \mathbb{E}_{\theta} \{\log q_{\theta}(Y) \, - \, \log \hat{q}_U(Y;X) \} \\
\, & = & \,  H(\hat{q}_U) \, - \, H(q_0)\, \\
\, & = & \,  H_d(a, \Sigma_S, \mathcal{L}_3) \, - \, H_d(a, \Sigma_Y, \mathcal{L}_2)\,,
\end{eqnarray*}
by the independence of $X$ and $Y$, the constancy of location family density $q_{\theta}$,  and since $Y-X|\theta \sim MMN_d(0,a, \Sigma_S, \mathcal{L}_3)$.  The result then follows from Lemma \ref{lemmaentropymmn}. \qed

The particular case with $\Sigma_X=\sigma^2_X I_d$ and $\Sigma_Y=\sigma^2_Y I_d$ follows directly from (\ref{riskmre}) and yields
\begin{equation}
\label{minimaxrisk}
R_{KL}(\theta, \hat{q}_U) \, = \, H_1(\frac{\|a\|}{\sigma_S}, 1, \mathcal{L}_3) \, - \, H_1(\frac{\|a\|}{\sigma_Y}, 1, \mathcal{L}_2) + \frac{d}{2} \log \frac{\sigma^2_S}{\sigma^2_Y}\,,
\end{equation}

\subsection{Minimum risk predictive density: Examples}

Theorem 2.1 tells us that the minimum risk predictive density is given by $\hat{q}_U(\cdot;X) \sim \hbox{MMN}_d(X, a, \Sigma_X + \Sigma_Y, \mathcal{L}_3)$ with $\mathcal{L}_3$ the cdf of $V_2-V_1$.   The result is quite general and can be viewed as an extension of the multivariate normal case with $a=0$ and $\hat{q}_U(\cdot;X) \sim \hbox{N}_d(X, \Sigma_X + \Sigma_Y)$.   Here are some interesting examples.   When continuous, the mixing distributions can be taken to have a scale parameter equal to one without loss of generality, since a multiple can be integrated into the shape vector $a$.  
\begin{enumerate}
\item[ (A)]   For the case of degenerate $V_2$ with $\mathbb{P}(V_2=v_2)=1$, i.e., when the distribution of $Y|\theta$ is normal with $Y \sim N_d(\theta+av_2, \Sigma_Y)$, the minimum risk equivariant predictive density reduces to $\hat{q}_U(\cdot;X) \sim MMN_d(X+av_2, -a , \Sigma_X+\Sigma_Y, \mathcal{L}_1)$.

\item[ (B)]  For the case of degenerate $V_1$ with $\mathbb{P}(V_1=v_1)=1$, i.e., when the distribution of $X$ is normal with $X|\theta \sim N_d(\theta+av_1, \Sigma_Y)$, the minimum risk equivariant predictive density reduces to $\hat{q}_U(\cdot;X) \sim MMN_d(X-av_1, a , \Sigma_X+\Sigma_Y, \mathcal{L}_2)$.

\item[ (C)]   We consider in this example $V_1, V_2$ i.i.d. exponentially distributed with densities $f(t)=e^{-t} \, \mathbb{I}_{(0,\infty)}(t)$, as well as $\Sigma_X=\sigma^2_X \, I_d$ and $\Sigma_Y=\sigma^2_Y \, I_d$.  
Here the distribution of $V_3$ is Laplace or double-exponential with density $\frac{1}{2} \, e^{-|v_3|}$ on $\mathbb{R}$.   Therefore, from Theorem 2.1, we have
\begin{eqnarray*}
\hat{q}_U(y;x)\, & = &\, \int_{\mathbb{R}} \frac{1}{2} \, e^{-|v_3|} \frac{1}{\sigma_S^d} \, \phi_d(\frac{y-x-av_3}{\sigma_S})  \, dv_3\,, \\
& = & \phi_d \left( y-x;\sigma^2_S I_d\right) \,  \int_{\mathbb{R}_+} 
e^{- (v_3^2 \frac{\|a\|^2}{2 \sigma^2_S} + v_3)} \;  \cosh \left( v_3 (\frac{(y-x)^{\top} a}{\sigma^2_S}\right) \, dv_3
\end{eqnarray*}
with $\sigma_S = (\sigma^2_X + \sigma^2_Y)^{1/2}$.  By making use of Lemma \ref{lemmaJ} in the Appendix with $A= \frac{\|a\|^2}{\sigma^2_S}$, $B= -1 \, \pm \,  \frac{(y-x)^{\top} a}{\sigma^2_S}$, and $c=0$, we obtain (for $a \neq 0$)
\begin{equation}
\nonumber
\begin{split}
\hat{q}_U(y;x)\,   =  & \, \sqrt{\frac{\pi \sigma^2_S}{2 \|a\|^2}} \, \phi_d(y-x;\sigma^2_S I_d) \, e^{\frac{\sigma^2_S}{2 \|a\|^2} \, + \, \frac{\{(y-x)^\top a\}^2}{2 \sigma^2_S \|a\|^2}}\\  \, & \times \, \left[\left\lbrace e^{- \frac{(y-x)^{\top} a}{\|a\|^2}} \, 
  \Phi\left( \frac{\sigma_S}{ \|a\|} (\frac{(y-x)^{\top} a}{\sigma^2_S} -1 ) \right) \right\rbrace + \left\lbrace e^{ \frac{(y-x)^{\top} a}{\|a\|^2}} \, 
\Phi\left(- \frac{\sigma_S}{ \|a\|} (\frac{(y-x)^{\top} a}{\sigma^2_S} +1 ) \right)  \right\rbrace \right] \,.
\end{split}
\end{equation}

\item[ (D)]   Consider $V_1, V_2$ i.i.d. truncated normal distributed $\hbox{TN}(0,1) $ (or equivalently as $\sqrt{\chi^2_1}$) for which $X$ and $Y$ are i.i.d. as multivariate skew-normal as in (\ref{skewnormaldensity}).   A straightforward calculation yields the density 
\begin{equation}
\nonumber g_{V_3}(t) \, = \, 2 \sqrt{2}  \; \phi(\frac{t}{\sqrt{2}}) \; \Phi(-\frac{|t|}{\sqrt{2}}) \,\, \mathbb{I}_{\mathbb{R}}(t)\,,
\end{equation} 
for $V_3 =^d V_1-V_2$.   It follows from Theorem 2.1, for $\Sigma_X=\sigma^2_X \, I_d$ and $\Sigma_Y=\sigma^2_Y \, I_d$, denoting $\sigma_S = (\sigma^2_X + \sigma^2_Y)^{1/2}$,  that
\begin{eqnarray}
\nonumber
\hat{q}_U(y;x)\, & = &\, \int_{\mathbb{R}} \, 2 \sqrt{2}  \; \phi(\frac{t}{\sqrt{2}}) \; \Phi(-\frac{t}{\sqrt{2}}) \, \phi_d(y-x-at; \sigma_S^2 I_d)  \, dt\,, \\
\label{mreTN}
\nonumber
&=&   \frac{2}{\sqrt{\pi}} \, \phi_d(y-x;  \sigma_S^2 I_d) \, 
\int_{\mathbb{R}_+} \Phi(-\frac{t}{\sqrt{2}}) \, e^{- \frac{t^2}{2} (\frac{1}{2} + \frac{a^{\top}a}{\sigma_S^2})} \left\lbrace e^{\frac{(y-x)^{\top}a t}{\sigma^2_S}} + e^{- \frac{(y-x)^{\top}a t}{\sigma^2_S}}  \right\rbrace   dt \,.
\end{eqnarray}
%The above can be evaluated with Lemma \ref{lemmaJ} with $c=-\frac{\sqrt{2}}{2}$, $A=\frac{2 \sigma_S^2}{\sigma_S^2+2a^{\top}a}$, and $B= \pm \frac{(y-x)^{\top}a }{\sigma^2_S}$, and collecting terms yielding
%\begin{equation}
%\hat{q}_U(y;x) \, = \, \phi_d(y-x;  \sigma_S^2 I_d) \, 
%\end{equation}

Now, by making use of Lemma \ref{lemmaJ} with $c=-\frac{\sqrt{2}}{2}$, $A=\frac{2 \sigma_S^2}{\sigma_S^2+2a^{\top}a}$, and $B= \pm \frac{(y-x)^{\top}a }{\sigma^2_S}$, collecting terms, and setting $f_k= \sqrt{\sigma^2_S + k a^{\top} a}$, we obtain the minimum risk equivariant predictive density
\begin{eqnarray*}
\nonumber
\hat{q}_U(y;x)\, & = & \,  \frac{4 \sigma_S}{f_1} 
\, \phi_d\left( y-x; \sigma^2_S (I_d + \frac{ a a^{\top}}{f_1^2})\right) \,  \\ & \! \times & \! \! \! \left\lbrace\Phi_2\left(- \frac{(y-x)^{\top} a}{f_1 f_2}, \frac{\sqrt{2}(y-x)^{\top} a}{\sigma_S f_2}; \frac{-\sigma_S}{\sqrt{2} \, f_2}\right) \, + \, 
\Phi_2\left(\frac{(y-x)^{\top} a}{f_1 f_2}, - \frac{\sqrt{2}(y-x)^{\top} a}{\sigma_S f_2}; \frac{-\sigma_S}{\sqrt{2} f_2}\right) \right\rbrace \,,
\end{eqnarray*}
where $\Phi_2(z_1, z_2; \rho)$ the cdf evaluated at $z_1, z_2 \in \mathbb{R}$ of a bivariate normal distributions with means equal to $0$, variances equal to $1$ and covariance equal to $\rho$.  In the evaluation above, we made use of the identities $(I-\frac{aa^{\top}}{f_2^2})^{-1} \, = \, I + \frac{aa^{\top}}{f_1^2}$ and $|I + \frac{aa^{\top}}{f_1^2}| \,= \, 1 +\frac{a^{\top}a}{f_1^2}\,$, which is a special case of the Sherman-Morrison formula for the matrix inversion of $A+b_1 b_2^{\top}$ with $A$ being a square matrix and $b_1$ and $b_2$ vectors of the same dimension.

\end{enumerate}

\section{Bayes posterior analysis and predictive densities}

\label{section on bayes predictive densities}
In this section, we expand on and document representations for Bayesian posterior and predictive densities for mean mixture of normal distributions.  

\subsection{Posterior densities}

Bayesian posterior analysis of MMN models relate to the general form 
\begin{equation}
\label{modelLM}  X|K,\theta \sim f_{\theta, K}\,,\, K \sim g\,,\, \hbox{ and } 
\theta \sim \pi\,,
\end{equation}
with observable $X \in \mathbb{R}^d$, density $g$ of $K$ free of $\theta$, and $\pi$ prior density for $\theta \in \mathbb{R}^d$.  Such a set-up leads to the following intermediate result, taken from \cite{lm2021}.

\begin{lemma}
\label{typeII}
For model (\ref{modelLM}), the posterior distribution of $U =^d \theta|x $ admits the representation
\begin{equation}
\label{typeIIposterior}   U|K' \sim \pi_{k',x} \hbox{ with } K' \sim g_{\pi,x}\,,
\end{equation}
$\pi_{k',x}$ being the posterior density of $\theta$ as if $K=k'$ had been observed, and $g_{\pi,x}(k') \propto g(k') \, m_{\pi,k'}(x) $ with $ m_{\pi,k'}$ being the marginal density of $X$ as if $K=k'$ had been observed.
\end{lemma}

We now apply the above to MMN distributions as in Definition \ref{defmmn}. 

\begin{example}
We apply Lemma \ref{typeII} to $X|\theta \sim MMN_d(\theta,a, \Sigma , \mathcal{L})$ and the prior $\theta \sim N_d(\mu, \Delta)$ with $\Sigma, \Delta>0$.  
The above fits into model (\ref{modelLM}) with $g$ taken to be the density of the mixing parameter $K =V \sim \mathcal{L}$, and $f_{\theta, k}$ the $N_d(\theta + k a, \Sigma)$ density.     Conditional on $K=k'$, standard Bayesian analysis for the normal model tells us that
\begin{equation}
  \theta|k',x \sim N_d \left( (I-P)x + P \mu - k'a, (I-P) \, \Sigma \right) \,, \, \hbox{ and } X|k' \sim N_d(\mu + k' a, \Sigma + \Delta)\,,
  \end{equation}
with $P=\Sigma \,(\Sigma+\Delta)^{-1}$, which yields the densities $\pi_{k',x}$ and $m_{\pi,k'}$ of Lemma \ref{typeII}.    Then from Lemma \ref{typeII}, we infer that
\begin{equation}
\label{posteriortheta}  
\theta|x \sim MMN_d \left((I-P)x + P \mu, \, a^* = - a \, , \, (I-P) \, \Sigma \,,  \mathcal{L}^*   \right)\,,
\end{equation}
where the distribution $\mathcal{L}^*$ has density
\begin{equation}
\label{posteriormixingdensity}
g_{\pi,x}(k') \, \propto   g(k') \, e^{-\frac{A}{2} k'^2 + Bk'}\,,\, \hbox{ with }  A= a^{\top} (\Sigma + \Delta)^{-1} a\, \hbox{ and } B= (x-\mu)^{\top} (\Sigma + \Delta)^{-1} a \,.
\end{equation}
Furthermore, it follows immediately that 
\begin{equation} 
\label{posteriorexpectation}
 \mathbb{E}(\theta|x) \, = \, (I-P)x + P \mu \, - \, P \, a \,  \mathbb{E}(K')\,, \, \hbox{ with } K' \sim g_{\pi,x}\,.
\end{equation}
\end{example}

\begin{remark}
For the improper prior density $\pi(\theta)\,=\,1$, one obtains $\theta|x \sim MMN_d(x,-a, \Sigma, \mathcal{L})$ by a direct calculation.  It can also be inferred from the above Example with $\Delta=\tau^2 I_d$ and $\tau^2 \to \infty$. 
\end{remark}

\begin{example}
It is interesting to further study the above posterior distributions for the particular cases where the mixing density (i.e., $V$ or $K$) of the MMN model is of the form  
\begin{equation}
\label{g}
g(k) \propto e^{-c_1 k^2/2 - c_2 k} \; \mathbb{I}_{(0,\infty)}(k), 
\end{equation}
with $c_1>0, c_2 \in \mathbb{R}$ or $c_1=0, c_2>0$.  Several of these distributions were presented in Example \ref{examplemixing}, but we recall that the cases $c_1>0$ for instance, which correspond to truncated normal distributions on $(0,\infty)$, lead to skew-normal densities (\ref{skewnormaldensity}) for $c_2=0$.
In the following, denote $\hbox{TN}\left( a,b; (0,\infty)\right)$ as a truncated normal distribution on $(0,\infty)$ with shape parameter $a \in \mathbb{R}$, scale parameter $b>0$, density $\frac{1}{b} \, \frac{\phi((y-a)/b)}{\Phi(a/b)} \, \mathbb{I}_{(0,\infty)}(y)$, and expectation $a \, + \, b R(a/b)$, with the reverse Mill's ratio $R(\cdot)$. 

Now, it is easily seen for cases where $K \sim g$ as in (\ref{g}) that
\begin{eqnarray*}
\label{gpi}
g_{\pi,x}(k') \, &\propto & e^{-(c_1+A) k'^2/2 \, + \, (B-c_2) k'} \, \mathbb{I}_{(0,\infty)}(k') \\
\,   &\propto &      \phi\left( \sqrt{A+c_1} \, k' \, - \, \frac{(B-c_2)}{\sqrt{A+c_1}} \right)     \mathbb{I}_{(0,\infty)}(k')\,,
\end{eqnarray*}
which is the density of a $\hbox{TN}\left(\frac{B-c_2}{A+c_1}, \frac{1}{\sqrt{A+c_1}}; (0,\infty)\right)$ distribution.   Hence, the above, which yields the density associated with $\mathcal{L}$, provides a complete description of the posterior distribution in (\ref{posteriortheta}) for all considered cases of mixing density (\ref{g}).  Analogously, the corresponding expectation  $\mathbb{E}(K') \, = \, \frac{B-c_2}{A+c_1} \, + \frac{1}{\sqrt{A+c_1}} \, R(\frac{B-c_2}{\sqrt{A+c_1}})$ provides an explicit expression for the posterior expectation $\mathbb{E}(\theta|x)$ in (\ref{posteriorexpectation}).

\end{example}
 
\subsection{Predictive densities} 

We now continue the above posterior analysis by focussing on the Bayes predictive density (i.e., the conditional density of $Y$ given $X=x$) for MMN distributions and a normally distributed prior for the unknown location parameter.    In doing so, the following extension come into play.

\begin{definition}
\label{defmmn2}
A random vector $Z \in \mathbb{R}^d$ is said to have a mean mixture of normal distribution  with two directions, denoted as $Z \sim MMN_d(\theta,a_1,a_2, \Sigma , \mathcal{L})$, if it admits the representation
\begin{align}
    \nonumber Z|V_1, V_2 &\sim N_d\left(\theta +  a_1W_1+a_2W_2\,, \Sigma\right) \hbox{ with }
    (W_1,W_2) \sim \mathcal{L},
    \label{mmnfamily}
\end{align}
where $\theta \in \mathbb{R}^d $ is a location parameter, $a_1, a_2 \in \mathbb{R}^d$ 
are known perturbation vectors, $\Sigma$ is a known positive definite covariance matrix, and $W_1,W_2$ are scalar random variable with joint cdf $\mathcal{L}$.
\end{definition}

We make use of the following intermediate result provided in \cite{lm2021} and applicable to mixture models of the form:
\begin{equation}
\label{*}
X|K,\theta \sim f_{\theta, K} \hbox{ with } K \sim g\,; Y|J,\theta \sim f_{\theta, J} \hbox{ with } J \sim h, \hbox{ and } \theta \sim \pi.  
\end{equation}
In the above set-up, $X \in \mathbb{R}^d$ is observable, the mixing variables $K$ and $J$ are independently distributed with distributions free of $\theta$, the variables $X$ and $Y$ are conditionally independent on $\theta$, and $\pi$ is a prior density for $\theta \in \mathbb{R}^d$ with respect to a $\sigma-$finite measure $\nu$.

\begin{lemma}
\label{typeIIpredictive}
For model (\ref{*}), setting $\pi_{k',x}$ and $g_{\pi,x}$ as in Lemma \ref{typeII}, the Bayes predictive density of $Y$ admits the mixture representation
\begin{equation}
\nonumber  Y|J',K' \sim q_{\pi}(\cdot|J',K'), \hbox{ with } J'  \sim h, K' \sim g_{\pi,x} \hbox{ independent },
\end{equation}
and $q_{\pi}(y|j',k') \,=\, \int_{\mathbb{R}^d} q_{\theta, j'}(y) \, \pi_{k',x}(\theta) \, d\nu(\theta)$, which can be interpreted as the Bayes predictive density for $Y$ as if $Y \sim q_{\theta, j'}$ and $K=k'$ had been observed.  
\end{lemma}

Applied to mean mixture of multivariate normal distributions with a normal distributed prior, we obtain the following presented as a theorem.
  
\begin{theorem}
\label{theorem-predictive}
\begin{enumerate}
\item[ {\bf (a)}]
For $X|\theta \sim MMN_d(\theta, a_X, \sigma_X^2 I_d, \mathcal{L}_1)$ and 
$Y|\theta \sim MMN_d(\theta, a_Y, \sigma_Y^2 I_d, \mathcal{L}_2)$ independent with prior $\theta \sim N_d(\mu, \tau^2 I_d)$, the Bayes predictive density for $Y$ is that of 
a 
$$MMN_d\left(\omega x + (1-\omega) \mu, - \omega a_X, a_Y, (\omega \sigma_X^2 + \sigma_Y^2) I_d, \mathcal{L}\right)\,$$
distribution, with $\mathcal{L}$ the joint cdf of $(K',J')$ with independently distributed $K' \sim g_{\pi,x}$ as in (\ref{posteriormixingdensity}) and $J' \sim \mathcal{L}_2$,
with $\omega=\tau^2/(\tau^2+\sigma_X^2)$, 
$A= \|a_X\|^2/(\sigma_X^2+\tau^2)$, and $B= \{(x-\mu)^{\top} a_X\}/(\sigma_X^2+\tau^2))$.

\item[ {\bf (b)}]  Moreover, whenever $a_Y= c a_X$ for $a_X \neq 0$ and a fixed $c \in \mathbb{R}$, the above predictive distribution is $MMN_d \, \left( \omega x + (1-\omega) \mu, a_X, (\omega \sigma_X^2 + \sigma_Y^2) I_d,  \mathcal{L}_3\right)$, with  $\mathcal{L}_3$ the cdf of $c J' - \omega K'$, and $(J',K')$ distributed as above.  Finally, for $a_X=0$, i.e., for $X|\theta \sim N_d(\theta, \sigma_X^2 I_d)$, the predictive distribution is $MMN_d(\omega x + (1-\omega) \mu, a_Y, (\omega \sigma_X^2 + \sigma_Y^2) I_d , \mathcal{L}_2)$
\end{enumerate}
\end{theorem}
{\bf Proof.}   Part {\bf (b)} follows immediately from part {\bf (a)}.   For part {\bf (a)}, consider $X'=X - K' a_X$ and $Y'=Y-J' a_Y$.  The result then follows from Lemma \ref{typeIIpredictive} with the familiar predictive density estimation result:
\begin{equation}
\nonumber   Y'|J', K',X' \sim N_d\left(\omega X' + (1-\omega) \mu, (\omega \sigma_X^2 + \sigma_Y^2) I_d \right)  ,
\end{equation}
implying 
$$  q_{\pi}(\cdot|J',K') \, \sim \, N_d\left(\omega x + (1-\omega) \mu - \omega a_X K' + a_Y J', (\omega \sigma_X^2 + \sigma_Y^2) I_d \right)\,, $$
matching Definition \ref{defmmn2} with $(W_1, W_2) =^d (K',J')$.
\qed

\begin{remark}
We point out that the minimum risk predictive density matches the density in {\bf (b)} with $\tau^2=\infty$, i.e., $\omega=1$.
\end{remark}

\section{Dominance Results}  

In this section, we first provide KL risk improvements on the MRE predictive density  $\hat{q}_U$ for estimating the density of $Y|\theta \sim MMN_d(\theta, a, \sigma^2_Y I_d , \mathcal{L}_2)$ based on $X|\theta \sim MMN_d(\theta, a, \sigma^2_X I_d , \mathcal{L}_1)$ with $d \geq 4$.  Such improvements are necessarily minimax as a consequence of Theorem \ref{theoremmre}.
Our findings cover two types of improvements: {\bf (i)} plug-in type (Section 4.1), and {\bf (ii)} Bayesian improvements (Section 4.2).   Furthermore, we provide analogue results for certain type of restricted parameter spaces which are also applicable for $d=2,3$.  Examples will be provided in Section 5.

\noindent   The restriction to covariance matrices that are multiple of identity is justified by convenience and the fact that there is no loss of generality in doing so.

\begin{remark} 
\label{intrinsic}  Predictive density estimates are intrinsic by nature which implies that the developments of this section, presented for $\Sigma_X=\sigma^2_X I_d$ and $\Sigma_Y=\sigma^2_Y I_d$ in model (\ref{model}) with known $\sigma^2_X$ and $\sigma^2_Y$, apply as well for $\Sigma_Y=c \Sigma_X$ with known $\Sigma_X, \Sigma_Y$, and $c=\sigma^2_Y/\sigma^2_X$.   
Indeed, one can consider $X'=\Sigma_X^{-1/2} X$ for which $X|\theta \sim MMN_d(\Sigma_X^{-1/2} \theta, \Sigma_X^{-1/2} a, I_d, \mathcal{L}_1)$ to estimate the density of $Y'=\Sigma_X^{-1/2} Y$, for which $Y'|\theta \sim MMN_d(\Sigma_X^{-1/2} \theta, \Sigma_X^{-1/2} a, c I_d, \mathcal{L}_2)$.  In doing so, one produces a predictive density estimator $q_1(y')\,=\, \hat{q}(y';x'), y' \in \mathbb{R}^d$, for the density $q_{Y'}$ of $Y'$, which equates to $q_2(y)\,=\, \hat{q}(\Sigma_X^{-1/2}y; \Sigma_X^{-1/2}x) \, |\Sigma_X^{-1/2}|$; $y \in \mathbb{R}^d$; as a predictive density estimator of the density $q_Y$ of $Y$.   Moreover, the Kullback-Leibler $\rho(q_{Y'},q_1)$ and $\rho(q_{Y},q_2)$ are equal, i.e.
\begin{equation}
\nonumber  
\int_{\mathbb{R}^d} q_{Y'}(t) \, \log \frac{q_{Y'}(t)}{q_1(t)} \, dt \,= \,  \int_{\mathbb{R}^d} q_Y(t) \, \log \frac{q_Y(t)}{q_2(t)} \, dt\,
,
\end{equation}
as seen with the change of variables $t \to \Sigma_X^{-1/2} t$.

\end{remark}

\subsection{Plug-in type improvements}

In the normal case with $X|\theta \sim N_d\left(\theta, \sigma^2_X I_d \right)$ and $Y|\theta \sim N_d \left(\theta, \sigma^2_Y I_d \right)$ independently distributed, the MRE predictive density $\hat{q}_U(\cdot;X) \sim N_d\left( X, (\sigma^2_X+ \sigma^2_Y) I_d \right)$ is inadmissible for $d \geq 3$ and can be improved by plug-in type densities of the form $q_{\hat{\theta}}(\cdot;X) \sim N_d \left(\hat{\theta}(X), (\sigma^2_X+ \sigma^2_Y) I_d \right)$.  Indeed, the KL risk performance 
of $q_{\hat{\theta}}$ relates directly to the ``dual'' point estimation risk of $\hat{\theta}(X)$ for estimating $\theta$ under squared error loss $\|\hat{\theta}-\theta\|^2$, with $q_{\hat{\theta}}(\cdot;X)$ dominating $\hat{q}_U(\cdot;X)$ if and only if $\hat{\theta}(X)$ dominates $X$ (\cite{fmrs2011}).  For MMN distributions, such a duality does not deploy itself in the same way, but does so after transformation of $(X,Y)$ to a canonical form and through the intrinsic nature of predictive densities.    The following result exhibits this and is applicable to $d \geq 4$.

\begin{theorem}
\label{theoremdominanceplugins}
 Consider $X,Y$ distributed as in model (\ref{modelXY}) with $ a \neq 0, d \geq 4, \theta \in \mathbb{R}^d, \Sigma_X=\sigma^2_X I_d, \hbox{ and } \Sigma_Y=\sigma^2_Y I_d$,
 and the problem of obtaining a predictive density estimator $\hat{q}(y;X)$, $y \in \mathbb{R}^d$, for the density of $Y$.   Let $H = \begin{pmatrix} {h_1^{\top}}\\ {H_2}\end{pmatrix}$ be an $d \times d$ orthogonal matrix such that $h_1 = \frac{a}{\|a\|}$.  Define the densities $$q_1(\cdot;X) \sim MMN_1 \left( h_1^{\top} X, \|a\|, (\sigma^2_X + \sigma^2_Y), \mathcal{L}_3\right) \hbox{ and } q_{2, \hat{\zeta}_2}(\cdot;X) \sim N_{d-1}\left(\hat{\zeta}_2 (H_2 X), (\sigma^2_X + \sigma^2_Y) I_{d-1}\right)\,. $$    Then, the predictive density 
 $q_{H, \hat{\zeta}_2}(y;X) \, = \,  q_1(h_1^{\top}y;X) \,\times \, q_{2, \hat{\zeta}_2}(H_2y;X) $, $y \in \mathbb{R}^d$, dominates $\hat{q}_U$ under KL loss if and only if $\hat{\zeta}_2(Z_2)$ dominates $Z_2$ as an estimator of $\zeta_2 \in \mathbb{R}^{d-1}$ under squared error loss $\|\hat{\zeta}_2 - \zeta_{(2)}\|^2$ and for the model $Z_2|\zeta_2 \sim N_{d-1}\left(\zeta_{(2)}, \sigma^2_X \, I_{d-1} \right)$.
 \end{theorem}
 
{\bf Proof.}    Set 
\begin{equation}
\label{x'y'}
X'\,=\, HX = \begin{pmatrix} {X_1'}\\ {X_{(2)}'}\end{pmatrix}\,,\, Y'\,=\, HY = \begin{pmatrix} {Y_1'}\\ {Y_{(2)}'}\end{pmatrix}\,, \hbox{ and }   \zeta\,=\, H\theta = \begin{pmatrix} {\zeta_1}\\ {\zeta_{(2)}}\end{pmatrix}\,,
\end{equation}
with $X_1'\,=\, h_1^{\top} X$, $X_{(2)}'=H_2 X$, $Y_1'\,=\, h_1^{\top} Y$, $X_{(2)}'=H_2 X$, $\zeta_1\,=\, h_1^{\top} \theta$, and $\zeta_{(2)}=H_2\theta$.
From Lemma \ref{lemmacanonical},  we have that $X_1'$, $X_{(2)}'$, $Y_1'$, and $Y_{(2)}'$ are independently distributed with $X_1' \sim MMN_1\left(\zeta_1, \|a\|, \sigma^2_X,\mathcal{L}_1\right)$, $Y_1' \sim MMN_1\left(\zeta_1, \|a\|, \sigma^2_Y,\mathcal{L}_2\right)$, $X_{(2)}' \sim N_{d-1}(\zeta_{(2)}, \sigma^2_X I_{d-1})$, and  $Y_{(2)}' \sim N_{d-1}(\zeta_{(2)}, \sigma^2_Y I_{d-1})$.   \\

Now consider the class of predictive densities of the form
\begin{equation}
q_{\hat{\zeta_2}}(y';X') \, = \, q_1(y_1';X_1') \times q_{2,\hat{\zeta_2}}(y_2';X_2')\,, y'=(y_1',y_{(2)}') \in \mathbb{R}^d,
\end{equation}
for estimating the density of $Y'$.  As in Remark \ref{intrinsic}, the Kullback-Leibler risk performance of $q_{H, \hat{\zeta}_2}(\cdot;X)$ for estimating the density of $Y$ is equivalent to the Kullback-Leibler risk performance of $q_{\hat{\zeta_2}}(\cdot;X')$ for estimating the density of $Y'$.  Furthermore, observe that the MRE density estimator $\hat{q}_U$ equates to density $q_{\hat{\zeta}_{2,0}}(\cdot;X')$ with $\hat{\zeta}_{2,0}(Y_2')\,=Y_2'$.
It thus follows, with the independence of the components of $Y'$ and $X'$, Lemma \ref{lemmaKL}, and setting $Z_2=X_{(2)}$ that  
\begin{eqnarray}
     \nonumber R_{KL}(\theta, \hat{q}_U) -  R_{KL}(\theta, q_{H,\hat{\zeta}_2}) & = & R_{KL}(\theta, q_{\hat{\zeta}_{2,0}}) \, - R_{KL}(\theta, q_{\hat{\zeta}_{2}})  \\
   \nonumber  & = & \mathbb{E} \log \left(\frac{q_1(Y_1';X_1')}{q_1(Y_1';X_1')}  \right)  \, + \, 
 \mathbb{E} \log \left(\frac{q_{2,\hat{\zeta}_2}(Y_2';X_2')}{q_{2,\hat{\zeta}_{2,0}}(Y_2';X_2')}  \right) 
\\
     \label{pluginriskdifference} & = & \frac{1}{2(\sigma^2_X + \sigma^2_Y)} \left(\mathbb{E}\;||\hat{\zeta}_2(Z_2) - \zeta_2||^2 - \mathbb{E}\;||Z_2 - \zeta_2||^2\right) \,,
\end{eqnarray}
which yields the result.   
\qed

The above dominance finding is quite general with respect to the specifications of $a, \mathcal{L}_1$, and $\mathcal{L}_2$ of model (\ref{modelXY}).  Furthermore, observe by examining (\ref{pluginriskdifference}) that the risk difference depends on $\theta$ only through $\zeta_{(2)}=H_2 \theta$ and this for any choice of $H_2$.  More strikingly as seen with (\ref{pluginriskdifference}), the risk difference does not depend on the mixing distributions $\mathcal{L}_1$ and $\mathcal{L}_2$ and can be simply described by a quadratic risk difference of point estimators which arise in a $(d-1)$ variate normal distribution problem.  
An illustration of Theorem \ref{theoremdominanceplugins} will be presented in Section \ref{illustrations}.

\subsection{Bayesian improvements}

We now focus on Bayesian predictive densities that dominate $\hat{q}_U$.  
In doing so, we work with canonical forms as in Lemma \ref{lemmacanonical}, apply the partitioning argument of Lemma \ref{lemmaKL}, and take advantage of known results for prediction in $(d-1)$ multivariate normal models. 
We consider a class of improper priors on $\theta$ which is the product measure of a (improper) uniform density over the linear subspace spanned by $a$ and a second component of the prior ($\pi_0$) supported on the subspace orthogonal to $a$.  The measure of this nature splits resulting Bayes predictive densities into independent parts and leads to a decomposition the KL risk in two additive parts. Hence, the dominance result is obtained by dominating the part of the KL risk corresponding to the orthogonal space to $a$, where transformed variables are $N_{d-1}$ distributed and where we can capitalize on known results.   Namely, the superharmonicity of $\pi_0$, or its associated marginal density or its associated square root marginal density, will suffice for dominance and minimaxity.

\begin{theorem}
\label{theoremmaindominance}
Consider $X,Y$ distributed as in model (\ref{modelXY})
with $\Sigma_X=\sigma^2_X I_d$, $\Sigma_Y=\sigma^2_Y I_d$, and $d \geq 2$.
% and $\sigma^2_Y\,=\, c \sigma^2_X $ for some known $c>0$.  Let $W = \frac{\sigma^2_X Y + \sigma^2_Y X}{\sigma^2_X + \sigma^2_Y}$ and 
Let $H = \begin{pmatrix} {h_1^{\top}}\\ {H_2}\end{pmatrix}$ be an $d \times d$ orthogonal matrix such that $h_1 = \frac{a}{\|a\|}$.  Let $X', Y'$, and $\zeta$ be defined as in (\ref{x'y'}) 
and consider prior densities of the form
\begin{equation}
    \pi(\theta) = \pi_0\left( \zeta_{(2)} \right).
    \label{priormuchoice1}
\end{equation}
\begin{enumerate}
\item[ {\bf (a)}]  Then,  the Bayes predictive density for $Y$ is given by 
\begin{equation}
\label{qpi}
\hat{q}_{\pi}(y;X) \, = \, \hat{q}'_{\pi}(H y; X')\,, y \in  \mathbb{R}^d,
\end{equation}
with $\hat{q}'_{\pi}(\cdot;x')$ the Bayes predictive density for $Y'$ based on $X'$, given by 
\begin{equation}
\label{99}  \hat{q}'_{\pi}(y'; X') \, = \, \hat{q}_U(y_1';X_1') \times \, \hat{q}'_{\pi_0}(y_{(2)}';X'_{(2)}) \,,
\end{equation}
with: {\bf (i)} $\hat{q}_U(\cdot;X_1')$ the MRE density, given in Theorem \ref{theoremmre}, of $Y_1' \sim MMN_1(\zeta_1, \|a\|, \sigma^2_Y, \mathcal{L}_2)$ based on $X_1' \sim MMN_1(\zeta_1, \|a\|, \sigma^2_X, \mathcal{L}_1)$, and {\bf (ii)} $\hat{q}'_{\pi_0}(\cdot;X_2')$ the Bayes predictive density for $Y'_{(2)} \sim N_{d-1}(\zeta_{(2)}, \sigma^2_Y I_{d-1})$ based on $X'_{(2)} \sim N_{d-1}\left(\zeta_{(2)}, \sigma^2_X I_{d-1}\right)$ and for prior density $\pi_0(\zeta_{(2)})$ for $\zeta_{(2)}$;

\item[ {\bf (b)}]  If $d \geq 4$, then $\hat{q}_{\pi}$ given in (\ref{qpi}) dominates the MRE $\hat{q}_U$, and is therefore minimax, if and only if $\hat{q}'_{\pi_0}(\cdot;X_2')$ dominates the MRE density for $Y'_{(2)}$ based on $X'_{(2)}$ given by a $N_{d-1}(X'_{(2)}, (\sigma^2_X+ \sigma^2_Y) I_{d-1})$ density.
\end{enumerate}

\end{theorem}
{\bf Proof.}
\begin{enumerate}
    \item[ {\bf (a)}]Eq. (\ref{qpi}) follows from the transformation of variables under the orthogonal matrix $H$. Note that the distribution of the transformed variables is
\begin{eqnarray}
 \nonumber  \:\: X^{'} \sim MMN_d(\zeta, a_0, \sigma^2_X I_d, \mathcal{L}_1) \\
   \nonumber  \hbox{and } \; \; Y^{'}  \sim MMN_d(\zeta, a_0, \sigma^2_Y I_d, \mathcal{L}_2),
\end{eqnarray}
where $a_0 = \left(\frac{a}{\|a\|}, 0,\ldots , 0\right)^{\top}$.
The prior of the form (\ref{priormuchoice1}) induces an improper uniform measure on $\zeta_1$ and independent  $\pi_0(\zeta_{(2)})$ on $\zeta_{(2)}$. Along with the conditional independence of $Y^{'}_{1}$ and $Y^{'}_{(2)}$ given $\zeta$, we get the Bayes predictive density as (\ref{99}).
%\begin{equation}
%\nonumber  \hat{q}'_{\pi}(y'; x') \, = \, \hat{q}_U(y_1';X_1') \times \, \hat{q}'_{\pi_0}%(y_{(2)}';x'_{(2)}) \,.
%\end{equation}
    \item[ {\bf (b)}] Observe that the MRE density estimator $\hat{q}_{U}(\cdot;X)$ corresponds to $\pi_0(\theta) = 1$, i.e., the improper uniform density on $\zeta_{(2)} \in \mathbb{R}^{d-1}$. By virtue of Lemma \ref{lemmaKL}, the KL risk difference between $\hat{q}_{U}(\cdot;X)$ and $\hat{q}_{\pi}(\cdot;X)$ is then expressed as
    \begin{eqnarray*}
    \label{KLriskdiffeq}
     \nonumber R_{KL}(\theta, \hat{q}_U)  -  R_{KL}(\theta, \hat{q}_{\pi}) & = & \mathbb{E}\;\text{log}\;\hat{q}_{\pi}(Y;X) -  \mathbb{E}\;\text{log}\;\hat{q}_{U}(Y;X)\\
   \nonumber  & = & \mathbb{E}\;\text{log}\;\hat{q}'_{\pi_0}(Y_{(2)}';X'_{(2)}) - \mathbb{E}\;\text{log}\;\hat{q}'_{U}(Y_{(2)}';X'_{(2)}) \\
    & = & R_{KL}(\zeta_{(2)}, \hat{q}'_U)  -  R_{KL}(\zeta_{(2)}, \hat{q}'_{\pi_0}),
\end{eqnarray*}
and part $\bf (b)$ follows.   \qed
\end{enumerate}

\begin{remark}
\label{remarkindependence of L1 and L2}
Theorem \ref{theoremmaindominance}'s dominance finding in part {\bf (b)} is unified with respect to the model settings $a$, $\mathcal{L}_1$ and $\mathcal{L}_2$, as well as the dimension $d \geq 4, \sigma_X^2$, and $\sigma_Y^2$.  Furthermore, as seen in the lines of the proof, the difference in risks between the predictive densities $\hat{q}_U$ and $\hat{q}_{\pi}$: (i) does not depend on the mixing $\mathcal{L}_1$ and $\mathcal{L}_2$, and (ii) depends on $\theta$ only through $\zeta_{(2)}=H_2 \theta$.
\end{remark}

Starting with \cite{komaki2001}, continuing namely with \cite{glx2006}, several Bayesian predictive densities $\hat{q}'_{\pi_0}(\cdot;X_2')$ have been shown to satisfy the dominance condition in part {\bf (b)} of the above Theorem.  Such choices lead to dominating predictive densities of $\hat{q}_U$.   
In \cite{glx2006}, analogously to the quadratic risk estimation problem with multivariate normal observables (e.g., \cite{stein1981,fsw2018}), sufficient conditions for minimaxity are conveniently expressed in terms of the marginal density of $Z \sim N_{d-1}(\zeta_{(2)}, \sigma^{2}I_{d-1})$ associated with density $\pi_0$ and given by
\begin{equation}
\nonumber  m_{\pi_0}(z,\sigma^2) = \int_{\mathbb{R}^{d-1}} \phi_{d-1}(z - \zeta_{(2)}, \sigma^2I_{d-1}) \, \pi_0(\zeta_{(2)})  \, d\zeta_{(2)} \, . 
\end{equation}
The superharmonicity of either $\pi_0$, $m_{\pi_0}(z,\sigma^2)$ for $z \in \mathbb{R}^{d-1}$,  for various values of $\sigma^2$, or as well of $\sqrt{m_{\pi_0}(z,\sigma^2)} $, each lead to sufficient conditions for minimaxity.   We recall here that the superharmonicity of $h: \mathbb{R}^{d-1} \to \mathbb{R}$ holds whenever the Laplacian $\Delta^2 h(t) = \sum_{i=1}^{d-1}\frac{\partial^2 h(t)}{\partial t_i^2}$ exists with $\Delta^2 h(t) \leq 0$ for $t \in \mathbb{R}^{d-1}$.

\begin{corollary}
\label{superharmonic conditions}
Consider the prediction context of Theorem \ref{theoremmaindominance} and a prior density $\pi_0$ as in (\ref{priormuchoice1}) other than the uniform density.  Suppose that $m_{\pi_0}(z, \sigma^2_X)$ is finite for all $z  \in   \mathbb{R}^{d-1}$ and that $d \geq 4$.  Then, the following conditions are each sufficient for $\hat{q}_{\pi}(\cdot;X)$ given in (\ref{qpi}) with prior density as in (\ref{priormuchoice1}) to dominate the MRE density $\hat{q}_U$: 
 \begin{enumerate}
     \item[{\bf (i)}]   $\Delta^2 \, m_{\pi_0}(z, \sigma^2) \leq 0$, $z \in \mathbb{R}^{d-1}$, for $\frac{\sigma^2_X\sigma^2_Y}{\sigma^2_X+\sigma^2_Y} < \sigma^2 < \sigma^2_X\,$, with strict inequality on a set of positive Lebesgue measure on $\mathbb{R}^{d-1}$ for at least one $\sigma^2$;
     \item[{\bf (ii)}]  $\Delta^2 \sqrt{m_{\pi_0}(z, \sigma^2)} \leq 0$,  $z \in \mathbb{R}^{d-1}$,  for  $\frac{\sigma^2_X\sigma^2_y}{\sigma^2_X+\sigma^2_Y} < \sigma^2 < \sigma^2_X\,$, with strict inequality on a set of positive Lebesgue measure on $\mathbb{R}^{d-1}$ for at least one $\sigma^2$;
     \item[{\bf (iii)}] The prior $\pi_0$ is such that $\Delta^2\pi_0(\zeta_{(2)}) \leq 0$ a.e. 
\end{enumerate}
      \end{corollary}
{\bf Proof.}  The results follow from part {(b)} of Theorem \ref{theoremmaindominance} and Theorem 1 - Corollary 2 in \cite{glx2006}.  \qed

Choices of the prior density $\pi_0$ satisfying the conditions of Corollary \ref{superharmonic conditions} thus rest upon analyses for the normal case which are plentiful.   In particular, several examples of $\pi_0$, and the resulting predictive density $\hat{q}_{\pi_0}'$, are provided in \cite{glx2006}.  These provide explicit representations of minimax predictive densities $\hat{q}_{\pi}$ given in (\ref{qpi}).  A detailed example is presented in Section 5.

%\begin{remark} 
The orthogonality decomposition used in this Section leads to a further interesting representation which generalizes the one obtained in the multivariate normal case, and for which we now expand upon.    For the multivariate normal case, referring to Theorem \ref{theoremmaindominance}'s decomposition, with  $X_{(2)}' \sim N_{d-1}(\zeta_{(2)}, \sigma^2_X I_{d-1})$ independent of $Y_{(2)}' \sim N_{d-1}(\zeta_{(2)}, \sigma^2_Y \, I_{d-1})$, a well-known representation of the Bayes predictive density associated with prior density $\pi_0$ for $\zeta_{(2)}$, given by \cite{glx2006}, is
\begin{equation}
\label{glxrepresentation}
    \hat{q}_{\pi_0}'(y_{(2)}';x_{(2)}')  = \hat{q}_{U}'(y_{(2)}';x_{(2)}') \times \frac{m_{\pi_0}(w_{(2)}'; \sigma^2_W)}{m_{\pi_0}( x_{(2)}',\sigma^2_X)} \,,
\end{equation}
with $w_{(2)}' = \frac{\sigma^2_X y_{(2)}' + \sigma^2_Y x_{(2)}'}{\sigma^2_X + \sigma^2_Y}$ and $\sigma^2_W = \frac{\sigma^2_X\sigma^2_Y}{\sigma^2_X+\sigma^2_Y},$ and where 
$\hat{q}_{U}'(\cdot;X_{(2)}')$ is the MRE predictive density of the density of $Y_{(2)}'$ based on $X_{(2)}'$, and given by a $N_{d-1}(x_{(2)}', (\sigma_X^2 + \sigma_Y^2) \, I_{d-1})$ density.

For the MMN case, we now have the following.
  
\begin{lemma}
\label{eqbayesrepresentation}
For a prior $\pi_0$ and $H$ in Theorem \ref{theoremmaindominance}, 
the corresponding Bayes predictive density $\hat{q}_{\pi}$ admits the representation
\begin{equation}
\label{glxextension}
    \hat{q}_{\pi}(y;x)  = \hat{q}_{U}(y;x) \times \frac{m_{\pi_0}(H_2 \, w,\sigma^2_W)}{m_{\pi_0}(H_2 \, x,\sigma^2_X)} \,,
\end{equation}
with $w = \frac{\sigma^2_X y + \sigma^2_Y x}{\sigma^2_X + \sigma^2_Y}$.
% and $\sigma^2_W = \frac{\sigma^2_X\sigma^2_Y}{\sigma^2_X+\sigma^2_Y}.$
\end{lemma}
{\bf Proof.}  Using the set-up of Theorem \ref{theoremmaindominance}, and expressions (\ref{qpi}) and (\ref{99}), the MRE predictive density is obtained as
$$ \hat{q}_U(y;X) \, = \, \hat{q}_U(y_1',X_1') \times  \hat{q}_{U}'(y_{(2)}';X_{(2)}')\,, y \in \mathbb{R}^d.$$   Therefore, from (\ref{qpi}) and (\ref{99}) again, as well as from \ref{glxrepresentation}, we obtain
$$  \hat{q}_{\pi}(y;X) \, = \, \hat{q}_U'(h_1^{\top} y; X_1') \times \hat{q}_U'(H_2 y; X_2') \times \frac{m_{\pi_0}(w_{(2)}'; \sigma^2_W)}{m_{\pi_0}( x_{(2)}',\sigma^2_X)} \,,  $$
which yields the result. \qed 

%\end{remark}
To conclude describing the dominance findings of this section and of Section 4.1,  we point out that the plug-in type improvements of Theorem \ref{theoremmaindominance} and the Bayesian dominance results of Theorem \ref{theoremmaindominance} and Corollary \ref{superharmonic conditions} are applicable regardless of the choice of the orthogonal completion $H_2$ of $H$, thus adding to choices of $\pi_0$ leading to minimaxity.   Furthermore, the above developments are unified and the findings are applicable for all MMN models (\ref{modelXY}) with $\Sigma_X=\sigma^2_X I_d$ and $\Sigma_Y=\sigma^2_Y I_d$, as well as for $\Sigma_Y=c \Sigma_X$ as justified in Remark \ref{intrinsic}. 

\begin{remark}
\label{remarkH2choice}
A particular appealing choice of $H_2$, which will be further explored below in Sections \ref{restrictedsection} and \ref{illustrations}, is the such that $H_2^{\top} H_2 \, = \, I_d - \frac{\; \; a a^{\top}}{a^{\top} a}$ in which case 
\begin{equation}
\|\zeta_{(2)}\|^2 \, = \, \theta^{\top} \left( I_d - \frac{\; \; a a^{\top}}{a^{\top} a} \right)  \theta\,,
\end{equation}
and spherically symmetric densities $\pi_{2}(\zeta_{2}) \, = \, g \left( \|\zeta_{2}\|^2 \right)$ lead to prior densities in (\ref{priormuchoice1}) of the form
\begin{equation}
\label{priorsphericallysymmetric}
\pi(\theta) \, = \, g \left\lbrace \theta^{\top} \left( I_d - \frac{\; \; a a^{\top}}{a^{\top} a} \right)  \theta \right\rbrace\, = \, g \left( \left\| \theta - \frac{a^{\top}\theta}{a^{\top}a} \, a  \right\| ^2  \right).
\end{equation}
\end{remark}
Such densities do not depend on $\|a\|$ and have contours given by hypersurfaces of cylinders with axis given by $a$ (or $h_1=\frac{a}{\|a\|}$).  Here is an example of three contours for $d=3$ and $a=(1,1,1)^{\top}$.

\begin{figure}[h]
         % \begin{subfigure}[b]{\textwidth}
         \centering
          \includegraphics[height=80mm,width=100mm]{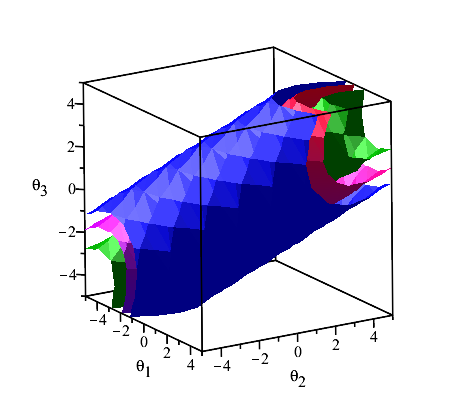}
\caption{Contours of $\pi(\theta)$ for $d=3$ and $a=(1,1,1)^{\top}$. }
\
\label{Fig99}
    % \end{subfigure}
   \end{figure}

\subsection{Restricted parameter spaces}
\label{restrictedsection}

Theorem \ref{theoremmaindominance}'s decomposition also leads to implications when there exists parametric restrictions on $\zeta_{(2)}\,=\, H \theta$.  Statistical models where parametric restrictions are present appear naturally in a great variety of contexts, and there is a large literature on related inferential problems, namely for a decision-theoretic approach (e.g., \cite{ms2004,vaneeden2006}).   Questions of predictive analysis under parametric restrictions are also of interest with findings obtained in \cite{ms2018,kmst2013,fmrs2011}.   Namely, for normal models, specifically model (\ref{modelXY}) with $a=0$, $\Sigma_X=\sigma^2_X I_d$,  $\Sigma_Y=\sigma^2_Y I_d$ with $\theta$ constrained to a convex set $C_0$ with non-empty interior, \cite{fmrs2011} showed that the Bayes predictive density associated with the uniform prior for $\theta$ on $C_0$ dominates the MRE predictive density under Kullback-Leibler loss.  The next results extends this finding to MMN models.

\begin{theorem}
\label{restrictedspacetheorem}
Consider $X,Y$ distributed as in model (\ref{modelXY})
with $\Sigma_X=\sigma^2_X I_d$, $\Sigma_Y=\sigma^2_Y I_d$, and $d \geq 2$. Let $C \subset \mathbb{R}^{d-1}$ be a convex set with non-empty interior, and let $\pi_C(\theta) = \pi_{0,U}(\zeta_{(2)}) = I_C(\zeta_{(2)})\,.$   Then $\hat{q}_{\pi_C}(\cdot;X) $ dominates $\hat{q}_{U}(\cdot;X) $ under KL risk and the restriction $\theta \in \{\theta \in \mathbb{R}^{d}\, : \, H_2\theta \in C \}$. 
\end{theorem}
{\bf Proof.}  As in Theorem \ref{theoremmaindominance} and the given proof, we infer that $\hat{q}_{\pi}$ given in (\ref{qpi}) with prior density $\pi(\theta) \, = \, \pi_{0}(\zeta_{(2)})$ for $\zeta_{(2)}=H_2 \theta$ dominates $\hat{q}_{U}$ if and only if $\hat{q}'_{\pi_{0}}(\cdot;X_2')$ dominates the MRE density for $Y_{(2)}' \sim N_{d-1}(\zeta_{2}, \sigma^2_Y I_{d-1})$.  \footnote{Said otherwise, part {\bf (b)} of Theorem \ref{theoremmaindominance} could have been stated for $d \geq 2$, but this would lead to knowingly vacuous conditions in the absence of a parametric restriction.}   
But, since this latter dominance holds precisely for density $\pi = \pi_C$ for the uniform density choice $\pi_0 = \pi_{0,U}$ as shown in \cite{fmrs2011}, the result follows.   
\qed

The setting of $C$ above is quite general and interesting examples includes balls and cones.   As earlier, the finding is unified and general to the MMN models.    
Here are two applications of Theorem \ref{restrictedspacetheorem}.

\begin{example}
\label{examplerestricted}
Suppose $d = 2$, $a = (1,1)^{\top}$, and the parametric restriction $\underline{c} \leq \theta_1 - \theta_2 \leq \bar{c}$, with $C=(\underline{c}, \bar{c})$ a strict subset of $\mathbb{R}$.  The MRE density  $\hat{q}_{U}(\cdot;X) $ is that of $MMN_2(X, a, (\sigma^2_X + \sigma^2_Y) I_2,  \mathcal{L}_3)$ distribution.  In the context of Theorem \ref{restrictedspacetheorem},  we have $\zeta_{(2)} \, = \, 
\frac{\theta_1 - \theta_2}{\sqrt{2}}$  and the prior density $\pi_C(\theta) = I_{C}(\theta_1-\theta_2)$.      Theorem \ref{theoremmaindominance} tells us
that the Bayes predictive density $\hat{q}_{\pi_C}$ dominates the MRE $\hat{q}_U$ with respect to KL loss and under the given parametric restriction. 
\footnote{In Example \ref{examplerestricted}, for the compact interval case say without loss of generality $\underline{c}=-m$ and $\bar{c}=m$, there exists a much larger class of dominating predictive densities obtained by replacing the uniform density for $\zeta_{(2)}$ by an even density $\pi_0$ supported on $(-m,m)$ that is increasing and logconcave on $(0,m)$.  This is established as in Theorem \ref{restrictedspacetheorem} and making use of Theorem 3.2 in \cite{fmrs2011}, which exploits a related point estimation finding in \cite{kubokawa2005}.}

\noindent  An explicit expression for $\hat{q}_{\pi_C}$  is available from Lemma \ref{eqbayesrepresentation} with $\pi_0$ the uniform $U(\frac{\underline{c}}{\sqrt2}, \frac{\bar{c}}{\sqrt2})$ density for $\zeta_{(2)}$.    As evaluated in \cite{kmst2013}, we obtain
\begin{eqnarray*}
\left( \frac{\sqrt{2}}{\bar{c} - \underline{c}} \right) m_{\pi_0}(z, \sigma^2) \, & = & \, \int_{\underline{c}/\sqrt{2}}^{\bar{c}/\sqrt{2}} \phi\left( z -\zeta_{(2)}, \sigma^2  \right) \, d \zeta_{(2)} \\
\, & = &   \Phi\left(\frac{z + \bar{c}/\sqrt{2}}{\sigma} \right)  \, - \, \Phi\left(\frac{z + \underline{c}/\sqrt{2}}{\sigma} \right)  \,,
\end{eqnarray*}
and (\ref{glxextension}) then yields
\begin{equation}
\nonumber
\hat{q}_{\pi_{C}}(y;x) \, = \, \hat{q}_U(y;x) \,\; \frac{\Phi\left(\frac{w + \bar{c}/\sqrt{2}}{\sigma_W} \right)  \, - \, \Phi\left(\frac{w + \underline{c}/\sqrt{2}}{\sigma_W} \right)}{\Phi\left(\frac{x + \bar{c}/\sqrt{2}}{\sigma_X} \right)  \, - \, \Phi\left(\frac{x + \underline{c}/\sqrt{2}}{\sigma_X} \right)} \,,  y \in \mathbb{R},
\end{equation}
\end{example}
with $w \, = \, \frac{ \sigma^2_X y \, + \, \sigma^2_Y x}{\sigma^2_X + \sigma^2_Y}$, $\sigma^2_W \, = \, \frac{\sigma^2_X\sigma^2_Y}{\sigma^2_X + \sigma^2_Y}$, and $\hat{q}_U$ the MRE density which is that of a $MMN_1(x, a, (\sigma^2_X + \sigma^2_Y),  \mathcal{L}_3)$ distribution.
\begin{example}
Theorem \ref{restrictedspacetheorem} applies for $\theta$ restricted to a cylinder of radius, say $m$, with the axis along the direction $a$, i.e., 
\begin{equation}
  \nonumber  C_m = \left\{\theta \in \mathbb{R}^d : \left\| \theta - \frac{a^{\top}\theta}{a^{\top}a} a  \right\| \leq m \right\};
\end{equation}
examples of which are drawn in Figure \ref{Fig99}. 
The dominating predictive density $\hat{q}_{\pi_{C_m}}$ is Bayes with respect to the uniform prior density on $C_m$, which corresponds to (\ref{priorsphericallysymmetric}) with $g(t)  \, = \, I_{(0,m)}(t)$.  An explicit expression for $\hat{q}_{\pi_{C_m}}$ can be derived from Lemma \ref{eqbayesrepresentation} with $\pi_0$ the uniform density on the ball $B_m = \{t \in \mathbb{R}^{d-1}: \|t \| \leq m \}$ and marginal density
\begin{eqnarray*}
m_{\pi_0}(z, \sigma^2) \, & = & \, \int_{B_m} \phi_{d-1}\left( z -\zeta_{(2)}, \sigma^2 I_{d-1} \right) \, d \zeta_{(2)} \\
\, & = &   F_{d-1, \frac{\|z\|^2}{\sigma^2}}(\frac{m^2}{\sigma^2})\,,
\end{eqnarray*}
with $F_{\nu,\lambda}$ the cdf of a $\chi_{\nu}^2(\lambda)$ distribution.    From (\ref{glxextension}), we thus obtain 
\begin{equation}
\nonumber
\hat{q}_{\pi_{C_m}}(y;x) \, = \, \hat{q}_U(y;x) \, \left(\frac{F_{d-1, \frac{\|H_2w|^2}{\sigma^2_W}}(\frac{m^2}{\sigma^2_W})}{F_{d-1, \frac{\|H_2x|^2}{\sigma^2_X}}(\frac{m^2}{\sigma^2_X})} \right), \,  y \in \mathbb{R}^d,
\end{equation}
with $\|H_2 t\|^2 \, = \, t^{\top} \left( I - \frac{\; \; a a^{\top}}{a^{\top} a} \right)  t \,$, for $t \in \mathbb{R}^d$, $w \, = \, \frac{ \sigma^2_X y \, + \, \sigma^2_Y x}{\sigma^2_X + \sigma^2_Y}$, $\sigma^2_W \, = \, \frac{\sigma^2_X\sigma^2_Y}{\sigma^2_X + \sigma^2_Y}$, and $\hat{q}_U$ the MRE density which is that of a $MMN_d(x, a, (\sigma^2_X + \sigma^2_Y) I_d , \mathcal{L}_3)$ distribution.
\end{example}

\section{Illustrations}
\label{illustrations}

We provide here illustrations of Theorems \ref{theoremdominanceplugins} and \ref{theoremmaindominance} accompanied by numerical comparisons and various observations.

\begin{figure}
    
     \begin{subfigure}{.45\textwidth}
          \includegraphics[height=80mm,width=90mm]{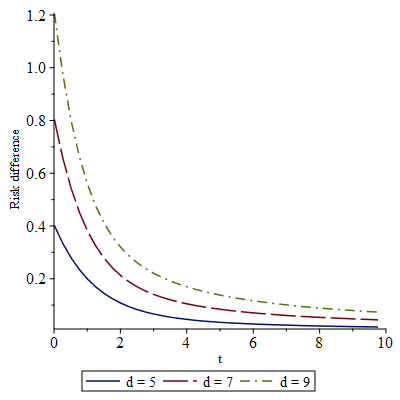}
\caption{ The KL Risk difference between $\hat{q}_U$ and $\hat{q}_{\pi_0}$ as a function of $t = \frac{\|\zeta_{(2)}\|^{2}}{d-1} = \frac{1}{d-1} \left\| \theta - \frac{a^{\top}\theta}{a^{\top}a} a  \right\| ^2  $, for $\sigma^2_Y=2, \sigma^2_X=1$. }
\
\label{KL risk diff Bayes}
     \end{subfigure}
     \hfill
     \begin{subfigure}{.45\textwidth}
          \includegraphics[height=80mm,width=90mm]{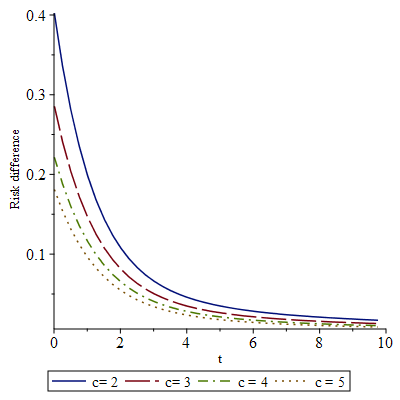}
\caption{The KL Risk difference between $\hat{q}_U$ and $\hat{q}_{\pi_0}$ as a function of $t = \frac{\|\zeta_{(2)}\|^{2}}{d-1} = \frac{1}{d-1} \left\| \theta - \frac{a^{\top}\theta}{a^{\top}a} a  \right\| ^2  $, for $d = 5$, $\sigma^2_X=1$, and  $c =  \frac{\sigma^2_Y}{\sigma^2_X}=2, 3, 4, 5$. }
\
\label{Riskdiffvaryingc}
     \end{subfigure} \\
     \vspace{15mm}
     \begin{subfigure}{.45\textwidth}
        \includegraphics[height=80mm,width=90mm]{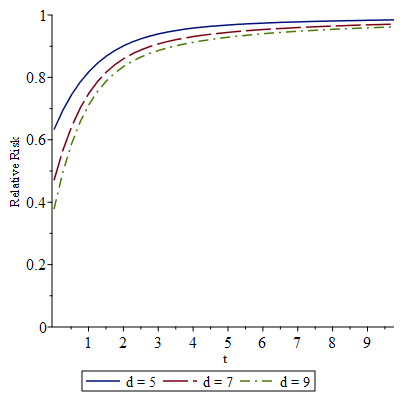}
\caption{Kullback-Leibler risk ratio between $\hat{q}_U$ and $\hat{q}_{\pi_0}$ as a function of $t = \frac{1}{d-1} \, ||  \theta - \bar{\theta} \boldsymbol{1}_d||  ^{2}  $, for $\sigma^2_Y=2, \sigma^2_X=1$. }
\
\label{FigrelativeriskBayes}
     \end{subfigure}
     \hfill
     \begin{subfigure}{.45\textwidth}
         \includegraphics[height=80mm,width=90mm]{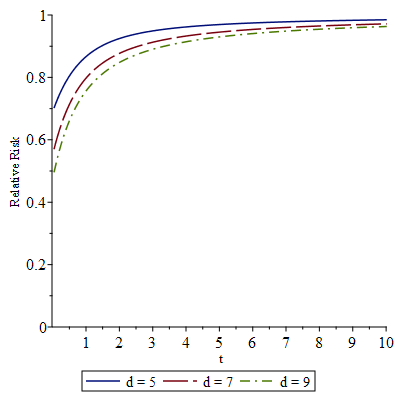}
\caption{ Kullback-Leibler risk ratio between $\hat{q}_U$ and $ q_1 \times q_{2,\hat{\zeta}_JS}$ as a function of $t = \frac{1}{d-1} \, ||  \theta - \bar{\theta} \boldsymbol{1}_d||  ^{2}  $, for $ \sigma^2_X = 1, \sigma^2_Y = 2, $ where $\hat{\zeta}_{JS}$ is James-Stein estimator.} 
\label{FigRelativeRiskJS}
     \end{subfigure}
       \caption{KL risk performance of the different predictive density estimators with the MRE}
\end{figure}

\begin{example} (A Bayesian minimax predictive density) 
\label{exampleBayes}
In the context of Theorem \ref{theoremmaindominance}, consider $H_2$ as in Remark \ref{remarkH2choice} combined with the harmonic prior density for $\zeta_{(2)} \in \mathbb{R}^{d-1}$ given by $\pi_0(\zeta_{(2)}) = \|\zeta_{(2)}\|^{-(d-3)}$ and which generates via (\ref{priorsphericallysymmetric}) an ``adjusted'' harmonic prior density on $\theta$ given by  
\begin{equation}
    \pi_H(\theta) =\left\| \theta - \frac{a^{\top}\theta}{a^{\top}a} a  \right\| ^{-(d-3)}\,.
    \label{harmonicpriorchoice}
\end{equation}

Thus, the prior density is the product measure on $\mathbb{R}^d$ with uniform prior on the linear subspace spanned by $a$ and the above harmonic measure on the $(d-1)-$dimensional chosen subspace orthogonal to $a$.   Since $\pi_0$ is superharmonic on $\mathbb{R}^{d-1}$ for $d \geq 4$, it follows from Corollary \ref{superharmonic conditions} that the Bayes predictive density $\hat{q}_{\pi_H}(\cdot;X)$ given in (\ref{qpi}), as well as in (\ref{qpiH}) below, dominates the MRE density $\hat{q}_U$ and is consequently minimax.

An explicit expression for $\hat{q}_{\pi_H}$ is available from Lemma \ref{glxextension} with marginal density 
\begin{align}
  \nonumber   m_{\pi_0}(z,\sigma^2) = \int_{\mathbb{R}^{d-1}} \phi_{d-1}(z - \zeta_{(2)}, \sigma^2I_{d-1}) \, \frac{1}{ ||\zeta_{(2)}||^{(d-3)}} \, d\zeta_{(2)} = \sigma^{3-d} \, \mathbb{E} \, T^{\frac{(3-d)}{2}},
\end{align} 
where $T \sim \chi^2_{d-1}\left( \frac{||z||^2}{\sigma^2}\right)$. In particular for odd $d \geq 5$, as shown in the Appendix, one may obtain
\begin{equation}
    \label{appendixformula}  m_{\pi_0}(z,\sigma^2) = \left(||z||^2\right)^{\frac{3-d}{2}}\left(1 - e^{-\frac{||z||^2}{2\sigma^2}}\sum\limits_{k = 0}^{\frac{d-5}{2}}\left(\frac{||z||^2}{2\sigma^2}\right)^{ k} \frac{1}{k!}\right) = r(||z||^2,\sigma^2)\, \hbox{ (say) },
\end{equation}
which relates to known results on the inverse moments of a chi-square variable with even degrees of freedom (e.g., \cite{bjy1984}), as well a closed form for an incomplete gamma function which intervenes in Komaki's \cite{komaki2001} representation of $m_{\pi_0}$.   From (\ref{glxextension}) and the above, we thus have 

\begin{equation}
\label{qpiH}
\hat{q}_{\pi_H}(y;x) \, =  \hat{q}_{U}(y;x) \,  \frac{r\left(\left\| w - \frac{a^{\top}w}{a^{\top}a} a  \right\| ^{2},\sigma^2_W\right)}{r\left(\left\| x - \frac{a^{\top}x}{a^{\top}a} a  \right\| ^{2},\sigma^2_X\right)} \,, y \in \mathbb{R}^d\,,
\end{equation}
 where $w$ and $\sigma^2_W$ are as given in  Lemma \ref{eqbayesrepresentation}.
%Thus, 
%\begin{eqnarray}
%   \nonumber \frac{\hat{q}_{\pi}(y;x)}{\hat{q}_{U}(y;x)} = \frac{r\left(\left\| w - \frac{a^{\top}w}{a^{\top}a} a  \right\| ^{2},\sigma^2_W\right)}{r\left(\left\| x - \frac{a^{\top}x}{a^{\top}a} a  \right\| ^{2},\sigma^2_X\right)},
%\end{eqnarray}

Risk differences between $\hat{q}_U$ and $\hat{q}_{\pi_H}$ are plotted in {\bf Figure
\ref{KL risk diff Bayes}} and {\bf Figure \ref{Riskdiffvaryingc}} as a function of $\|\zeta_{(2)}\|^2$, or equivalently as a function of 
\begin{eqnarray}
   \nonumber  t = \frac{\|\zeta_{(2)}\|^{2}}{d-1} = \frac{1}{d-1} \left\| \theta - \frac{a^{\top}\theta}{a^{\top}a} a  \right\| ^2 \, ,
\end{eqnarray}
i.e., in terms of the average squared component of $\zeta_{(2)}$.  The actual risks depend on the underlying mixing distributions $\mathcal{L}_1$ and $\mathcal{L}_2$, but not the risk differences as previously observed in Remark \ref{remarkindependence of L1 and L2}.   Observe as well that $t$ is independent of $\|a\|$ and only depends on the direction $a/\|a\|$.    Figure \ref{KL risk diff Bayes} has $\sigma^2_X = 1, \sigma^2_Y = 2 $ and varying $d$, while Figure \ref{Riskdiffvaryingc} has fixed $d=5, \sigma^2_X=1$ with $\sigma^2_Y=c \sigma^2_X$ and varying $c$.  
As seen with {\bf Figure \ref{KL risk diff Bayes}}, the improvement in KL risk vanishes at $t \to \infty$, but gains in prominence with increasing $d$, and with the proximity of $\theta$ to the linear subspace spanned by $a$.  As seen with {\bf Figure \ref{Riskdiffvaryingc}}, the KL risk difference loses in prominence with larger $c$ which is consistent with the fact that MRE density gains in reliability when the variance $\sigma^2_X$ of the observable decreases. 
 
Frequentist risk ratios between $\hat{q}_U$ and $\hat{q}_{\pi_H}$ are plotted in
{\bf Figure \ref{FigrelativeriskBayes}} for $\sigma^2_X = 1, \sigma^2_Y = 2 $ and varying $d$.  These ratios depend additionally on the mixing distributions $\mathcal{L}_1$ and $\mathcal{L}_2$ and they are set here with $\sqrt{\chi^2_1}$ mixing (Example 2.1 (B)), i.e., $X|\theta$ and $Y|\theta$ have skew-normal distributions with densities given in (\ref{skewnormaldensity}) and MRE density expanded upon in part (D) of Section 2.4.    We further set $a = \boldsymbol{1}_d = (1,\ldots,1)^{\top}$, in which case the harmonic prior density on $\theta$ in (\ref{harmonicpriorchoice}) reduces to $ \pi_0(\theta) = ||  \theta - \bar{\theta} \boldsymbol{1}_d||  ^{-(d-3)}$ with $\bar{\theta} = \frac{1}{d}\sum\limits_{i=1}^d \theta_i$. 
With the above settings, the constant (and minimax) risk of the MRE density can be computed from (\ref{minimaxrisk}).   For instance, we obtain $R(\theta, \hat{q}_U)  \approx 1.0954$ for $d=5$, $\approx 1.5187$ for $d=7$ and $\approx 1.9403$ for $d=9$.
These are close to linear with the term $\frac{d}{2} \log \frac{\sigma^2_S}{\sigma^2_Y} = \frac{d}{2} \log \frac{3}{2}$ ($\approx 1.0137$ for $d=5$, $\approx 1.4191$ for $d=7$ and $\approx 1.8246$ for $d=9$), representing the MRE risk for the normal case with $a=0$, being dominant in (\ref{minimaxrisk}).   As seen in {\bf Figure \ref{FigrelativeriskBayes}}, where the risk ratios are plotted with respect to $t = \frac{1}{d-1} \, ||  \theta - \bar{\theta} \boldsymbol{1}_d||  ^{2}  $, the gains increase in $d$ and with the closeness of the $\theta_i$'s to $\bar{\theta}$. 

%The risk expression reduces to
%\begin{align}
%   R_{KL}(\hat{q}_U,0) & = \nonumber = \mathbb{E}\log q(Y'_1) - \mathbb{E}\log %\hat{q}(Z'_1) - \frac{d}{2}\log \frac{\sigma^2_y}{\sigma^2_x  + \sigma^2_y}.
%\end{align}
%Here, $Y'_1 \sim SN_1(0,\sqrt{\frac{d}{\sigma^2_y}},1)$ and $Z'_1 \sim SN_1(0,%\sqrt{\frac{d}{\sigma^2_y + \sigma^2_x}},1)$.

\end{example}

\begin{example} (Plug-in type improved predictive density)   
In the context of Theorem \ref{theoremdominanceplugins}, consider plug-in type predictive densities $q_{H, \hat{\zeta}_2}(y;X)\, = \,  q_1(h_1^{\top}y;X) \,\times \, q_{2, \hat{\zeta}_2}(H_2y;X) $ with the choice of the James-Stein estimator $ \hat{\zeta}_2(Z_2) = \left(1 - \frac{(d-3) \sigma^2_X}{||Z_2||^2}\right) \,Z_{2}$ leading to the dominance of $q_{H, \hat{\zeta}_2}$ over $\hat{q}_U$ for $d \geq 4$.        
Both the dominating predictive density $q_{H, \hat{\zeta}_2}$ and the actual difference in risks do depend on the choice of $H_2$, but the KL risk difference, as given in (\ref{pluginriskdifference}) and mentioned at the end of Section 4.1, is independent of the underlying mixing distributions and will thus coincide with the corresponding difference stemming for $d-1$ dimensional normal models and which have appeared many times in the literature.   The difference in risks will be a function of $\zeta_{(2)}=H_2 \theta$ in general, and more precisely as a function of $\|\zeta_{(2)}\|^2$ in this case given that the James-Stein estimator is equivariant with respect to orthogonal transformations.
 
It is thus more interesting to look at the ratio of Kullback-Leibler risks and such ratios are  presented in {\bf Figure \ref{FigRelativeRiskJS}} with the same settings as in  Example \ref{exampleBayes}, i.e.,   multivariate skew-normal models with $\sqrt{\chi^2_1}$ mixing, $\sigma^2_X = 1,\sigma^2_Y = 2$, and $a = (1,\cdots,1)^T.$
Again here, the risk ratios are plotted with respect to $t = \frac{1}{d-1} \, ||  \theta - \bar{\theta} \boldsymbol{1}_d||  ^{2}  $, the gains increase in $d$ and with the closeness of the $\theta_i$'s to $\bar{\theta}$.
\end{example}

\newpage
\section*{Concluding remarks}

In this work, we have addressed the problem of determining efficient predictive densities under Kullback-Leibler frequentist risk for multivariate skew-normal distributions and, more generally, for mean mixtures of multivariate normal (MMN) distributions, and provided Bayesian and plug-in type predictive densities which dominate the MRE density, and are minimax in four dimensions or more.   In doing so, we have made use of a canonical transformation which leads to the decomposition of the Kullback-Leibler risk for the predictive densities being considered into two additive parts, one of which matching that of the MRE and minimax density, the other relating to a normal model and permitting improvement in view of shrinkage predictive  density estimation results for such models.     Further implications are provided for certain type of parametric restrictions.  In addition, motivated by the relative paucity of analytical representations for Bayesian posterior and predictive densities, we have contributed such explicit representations. 

This work represents, to the best of our knowledge, a first foray of the study of predictive density estimation for MMN distributions.   The findings are thus novel and they are also unified.  The canonical transformation technique may well find further applications in predictive analysis, such as for mean-variance mixture of normal distributions.    Extensions to other choices of loss (e.g., $\alpha$-divergence) and to unknown covariance structures would be most interesting to investigate as well.

\section*{Acknowledgements}
\'Eric Marchand's research is supported in part by the Natural Sciences and Engineering Research Council of Canada.   Pankaj Bhagwat is grateful to the ISM (Institut des sciences mathématiques) for financial support.   Thanks to Jean-Philippe Burelle for useful discussions on geometric representations related to prior density (\ref{priorsphericallysymmetric}).

\section*{Appendix}
%\begin{enumerate}
\begin{lemma}
\label{lemmaJ}
For all $B,c \in \mathbb{R}$, $A \in \mathbb{R}_+$, we have
\begin{equation} 
\int_0^{\infty} \Phi(c t) \, e^{- \frac{t^2}{2A} + Bt} \, dt \, = \, e^{\frac{AB^2}{2}} \, \sqrt{2\pi A} \,\, \Phi_2(\frac{c AB}{\sqrt{1+c^2 A}}, B \sqrt{A}; \frac{c \sqrt{A}}{\sqrt{1+c^2 A}}) \,. 
 \end{equation}
\end{lemma}    
{\bf Proof.}   We have  
\begin{eqnarray*}
e^{\frac{-AB^2}{2}} \, (2\pi A)^{-1/2} \, \int_0^{\infty} \Phi(c t) \, e^{- \frac{t^2}{2A} + Bt} \, dt \,
 \, & =  &  \int_0^{\infty} \Phi(c t) \, \frac{1}{\sqrt{A}}\,  \phi(\frac{t-AB}{\sqrt{A}}) \, dt\\
\, & = &   \mathbb{P}\left( U - cT \leq 0, -T\leq 0 \right)\, ,
\end{eqnarray*}
with $U,T$ independently distributed as $N(0,1)$ and $N(\theta_T=AB, \sigma^2_T=A)$, respectively. The result follows since
\begin{equation}
\nonumber
(U-cT, -T)^{\top} \sim N_2\left(  \left(\begin{array}{r} -cAB \\ AB
\end{array} \right) ,  \left[\begin{array}{rr}
1 +c^2 A& cA   \\
cA & A    
\end{array}  \right] \right)\,.  \;\qed
\end{equation}

{\bf Proof of (\ref{appendixformula}).}
%\end{enumerate}
With the standard representation $T |K  \sim \chi^2_{d-1 + 2K}$ with $K \sim Poisson\left( \frac{||z||^2}{2\sigma^2}\right)$, we have 
\begin{eqnarray*}
    \mathbb{E} \, T^{\,(3-d)/2} & = &\sum\limits_{k = 0}^{\infty} e^{-\frac{||z||^2}{2\sigma^2}}\frac{1}{k!}\left(\frac{||z||^2}{2\sigma^2}\right)^{k}\mathbf{E} \left(\chi^2_{d - 1 + 2k}\right)^{\frac{(3-d)}{2}} \\ 
%   & = & \sum\limits_{k = 0}^{\infty} e^{-\frac{||z||^2}{2\sigma^2}}\frac{1}{k!}\left(\frac{||z||^2}{2\sigma^2}\right)^{k} \frac{1}{2^{\frac{d-3}{2}}}\frac{\Gamma(k + 1)}{\Gamma(\frac{d-1}{2} + k)} \\
    & = & \frac{1}{2^{\frac{d-3}{2}}} e^{-\frac{||z||^2}{2\sigma^2}}\sum\limits_{k = 0}^{\infty} \left(\frac{||z||^2}{2\sigma^2}\right)^{k} \frac{1}{\Gamma(\frac{d-1}{2} + k)} \\
%        & = & \frac{1}{2^{\frac{d-3}{2}}} \, e^{-\frac{||z||^2}{2\sigma^2}} \left(\frac{||z||^2}{2\sigma^2}\right)^{-\frac{d-3}{2}} \sum\limits_{k = 0}^{\infty} \left(\frac{||z||^2}{2\sigma^2}\right)^{\frac{d-3}{2} + k} \frac{1}{\Gamma(\frac{d-1}{2} + k)} \\    
      & = & e^{-\frac{||z||^2}{2\sigma^2}} \left(\frac{||z||^2}{2\sigma^2}\right)^{-\frac{d-3}{2}} \sum\limits_{k = \frac{d-3}{2}}^{\infty}\left(\frac{||z||^2}{2\sigma^2}\right)^{ k} \frac{1}{k!}\,, 
\end{eqnarray*}
which yields (\ref{appendixformula}).  \qed

\end{document}